\documentclass{article}

\usepackage{PRIMEarxiv}
\usepackage{amsthm,amssymb}
\usepackage[utf8]{inputenc} % allow utf-8 input
\usepackage[T1]{fontenc}    % use 8-bit T1 fonts
\usepackage{hyperref}       % hyperlinks
\usepackage{url}            % simple URL typesetting
\usepackage{booktabs}       % professional-quality tables
\usepackage{amsfonts}       % blackboard math symbols
\usepackage{nicefrac}       % compact symbols for 1/2, etc.
\usepackage{microtype}      % microtypography
\usepackage{lipsum}
\usepackage{fancyhdr}       % header
\usepackage{graphicx}       % graphics
\graphicspath{{media/}}     % organize your images and other figures under media/ folder
\usepackage[fleqn]{amsmath}
\usepackage{paralist}
\usepackage[misc]{ifsym}
\usepackage{epsfig} 
\usepackage{epstopdf} 
\usepackage{tikz}
 \usepackage{pgfplots}
     \pgfplotsset{
        table/search path={TikZData},
    }
\pgfplotsset{compat=newest}

\newtheorem{theorem}{Theorem}[section]

\newtheorem{lemma}[theorem]{Lemma}

\newtheorem{example}[theorem]{Example}

\newtheorem{remark}[theorem]{Remark}

\newcommand{\R}{\mathbb{R}}
\newcommand{\de}{\mathrm{d}}
\newcommand{\DL}{\text{DL}}
\newcommand{\e}{\mathrm{e}}

\newenvironment{prf}[1][]{\noindent\textbf{Proof. #1}}{}

\usepackage{cases}   
\title{A spectral Galerkin exponential Euler time-stepping scheme for parabolic SPDEs on two-dimensional domains with a $\mathcal{C}^2$ boundary
}

\author{
  Julian Clausnitzer \\
  J\"{u}lich Supercomputing Centre, Forschungszentrum J\"{u}lich GmbH, 52425 J\"{u}lich \\
  {j.clausnitzer@fz-juelich.de} \\
   \AND
  Andreas Kleefeld \\
    J\"{u}lich Supercomputing Centre, Forschungszentrum J\"{u}lich GmbH, 52425 J\"{u}lich\\
  Faculty of Medical Engineering and Technomathematics, University of Applied Sciences Aachen, 52428 J\"{u}lich\\
  {a.kleefeld@fz-juelich.de} \\
}

\begin{document}
\maketitle

% Enter the first author's name and email address; email addresses are required for each author.
% Use footnote notations to indicate address and affiliations with commas between numbers if more than one address applies; see below for examples.

%%%%%%%%%%%%%%%%%%%%%%%%%%%%%%%%%%%%%%%%%%%%%%%%%%%%%%%
%             5. ABSTRACT
%%%%%%%%%%%%%%%%%%%%%%%%%%%%%%%%%%%%%%%%%%%%%%%%%%%%%%%

\begin{abstract}
We consider the numerical approximation of second-order semi-linear parabolic stochastic partial differential equations interpreted in the mild sense which we solve on general two-dimensional domains with a $\mathcal{C}^2$ boundary with homogeneous Dirichlet boundary conditions.
The equations are driven by Gaussian additive noise, and several Lipschitz-like conditions are imposed on the nonlinear function.
We discretize in space with a spectral Galerkin method and in time using an explicit Euler-like scheme. For irregular shapes, the necessary Dirichlet eigenvalues and eigenfunctions are obtained from a boundary integral equation method. This yields a nonlinear eigenvalue problem, which is discretized using a boundary element collocation method and is solved with the Beyn contour integral algorithm. We present an error analysis as well as numerical results on an exemplary asymmetric shape, and point out limitations of the approach.
\end{abstract}
%%%%%%%%%%%%%%%%%%%%%%%%%%%%%%%%%%%%%%%%%%%%%%%%%%%%%%
%                   6. BODY
%%%%%%%%%%%%%%%%%%%%%%%%%%%%%%%%%%%%%%%%%%%%%%%%%%%%%%
% Only the first word and proper nouns of section titles should be capitalized.
% The title of section 1:
 \section{Introduction.}
We consider the numerical approximation of semilinear parabolic stochastic partial differential equations (SPDEs) of the form
\begin{align*}
\de U(t) &= [\Delta U(t) + F(U(t))]\de t + \de W^Q(t), \\
U(0) &= u_0\in H, \quad t\in[0,T],
\end{align*}
where $T>0$, $H\subset L^2(D)$ is a Hilbert space of functions satisfying a homogeneous Dirichlet boundary condition on a domain $D\subset \R^2$ with a $\mathcal{C}^2$ boundary, $\Delta$ is the Laplace operator, $F$ is a nonlinear continuous operator and $W^Q$ is a cylindrical $Q$-Wiener process with covariance operator $Q$ (for more precise definitions, see \S \ref{setting}). We make the common assumption that the eigenfunctions of $Q$ and $\Delta$ coincide, which allows for a series expansion of $W^Q$ in terms of eigenfunctions of $\Delta$ (see Appendix \ref{appendixa}). \\
The numerical solution of SPDEs has attracted much attention, and many different methods for temporal as well as spatial discretization have been developed (for a comprehensive overview, see e.g. \cite[Section 1.2]{PusnikPhD}). However, due to the 
irregular behavior of SPDEs, strong regularity conditions have to be imposed on either the noise or the spatial domain.
Some publications focused on a one-dimensional spatial domain and potentially space-time white noise
\cite{Alabert,allennovoselzhang,blomker,brehier,daviegaines,gradinaru,winkel,KloedenRoyal,katsoulakis,walshcomeback}, while others 
have worked with trace-class or smooth noise and a potentially two-dimensional domain \cite[6.24]{chow}\cite{grecksch,hausenblas,jentzenpathwise,higherorder,lord2}.
We intend to use the latter approach by considering additive noise that has sufficiently quickly decaying modes, but in exchange being able to consider a much larger class of general two-dimensional domains $ D\subset \R^2$ with a $\mathcal{C}^2$ boundary. This broadening of the class of spatial domains has been pointed out in \cite[p. 260]{PusnikPhD} as a subject of future research. Also, while there has been some theoretical treatment of SPDEs on domains with a $\mathcal{C}^2$ boundary of dimension two or more with homogeneous Dirichlet boundary conditions \cite{bonaccorsi,dinunno,dozzi}, to our knowledge no solutions of SPDEs have been numerically approximated in such domains.  \\
Spectral methods have been widely used for the solution of semilinear parabolic SPDEs \cite{blomker,grecksch,winkel,shardlow} and are straightforward to implement in cases where the base functions are known explicitly. If this is not the case, the spectrum of the linear differential operator as well as its eigenfunctions have to be numerically approximated. This can be achieved using boundary element methods (BEM) resulting in a nonlinear eigenvalue problem. The advantage of this method is that any bounded domain with a $\mathcal{C}^2$ boundary can be considered, and the eigenfunction can be computed pointwise at any desired point in the interior of the domain, making it possible to store the functions on a uniform grid. Another possibility to compute the base functions are finite element methods (FEM), and for the Dirichlet Laplacian this has been done in \cite{femeigenvalue}. The difference to BEM is that instead of domains with a $\mathcal{C}^2$ boundary, FEM can deal with polygons with a Lipschitz boundary in $\R^2$. Since FEM requires the triangularization of the domain, if the functions are to be stored on a uniform grid, an additional error would be introduced by evaluating the functions on the triangularization at the required grid points. In this work, we focus on BEM and domains with $\mathcal{C}^2$ boundaries.\\
We give in \S\ref{setting} a precise mathematical framework, including the conditions placed on each component of the equations which can be solved. In \S\ref{section_Numericalscheme}, the discretization in space and time of the mild solutions is given as a Galerkin projection in space and the exponential Euler scheme in time as described by Jentzen and Kloeden \cite{KloedenRoyal}. The boundary integral equation method is recalled, including the necessary discretization steps turning the Helmholtz equation into a nonlinear eigenvalue problem. The algorithm developed by Beyn \cite{beyn}, which allows for the solution of nonlinear eigenvalue problems, is recalled, again also sketching the discretization steps necessary for the implementation. Lastly, some remarks on assembling the first elements of an orthonormal basis are made and it is shown how many boundary elements are needed to compute the eigenfunctions to a given accuracy. In \S\ref{section_erroranalysis}, we prove a result on the strong error of the spectral Galerkin-exponential Euler scheme with eigenvalues and eigenfunctions which are approximated to a given accuracy. In \S\ref{numericalexperiments}, we demonstrate the convergence on an asymmetric shape and present a few solutions for different nonlinear functions. The MATLAB program code, including the data and programs to generate the plots shown in the figures, is available at GitHub \href{https://github.com/JulianSPDE/2D-SPDE}{https://github.com/JulianSPDE/2D-SPDE}.
\section{Mathematical setting and assumptions.}\label{setting}
Let $D\subset \mathbb{R}^2$ be a bounded domain whose boundary is of class $\mathcal{C}^2$ and can be described by a parametrization $\gamma$ such that $\| \gamma'(t) \|_2>0$, $t\in[0,2\pi)$, $\gamma(0)=\gamma(2\pi)$ and $\gamma$ is nonintersecting, i.e. $\gamma(s)\neq \gamma(t)$ for distinct $s,t\in[0,2 \pi)$. Let $T>0$ be fixed. We consider SPDEs of the form
\begin{align}\label{spde}
\begin{split}
\de U(t) = [AU(t) + F(U(t))]\de t + \de W^Q(t),\quad
U(0) = u_0, \quad t\in[0,T]
\end{split}
\end{align}
where $(U(t))_{t\in [0,T]}$ is a stochastic process taking values in a Hilbert space $H$, \mbox{$u_0\in H$}, $A:D(A)\subset H\rightarrow H$ is a linear sectorial operator, i.e. it generates an analytic semigroup $(\e^{At})_{t\in[0,T]}$, and
$F: H\rightarrow H$ is a nonlinear operator. The process $(W^Q(t))_{t\in[0,T]}$ is a cylindrical $Q$-Wiener process taking values in $H$, where $Q:H\rightarrow H$ is a 
trace-class operator. We will consider the Laplacian $A = \Delta: D(\Delta)\subset H\rightarrow H$ with homogeneous Dirichlet boundary conditions on the domain $D$, so that
we have $H = L^2(D)$ and $D(A)=H_0^1(D)\cap H^2(D)$. It is also assumed that $\Delta$ and the covariance operator $Q$ commute, so that they have 
a complete set consisting of the same orthonormal eigenfunctions $e_i$, $i\in\mathbb{N}$. This assumption of a joint eigenbasis for $Q$ and $\Delta$ has been treated in \cite{blomkernoise} and we show in Appendix \ref{appendixa} that this assumption can be satisfied for more general two-dimensional shapes.\\
The SPDE \eqref{spde} can be rigorously understood as an integral equation involving a stochastic $\de W^Q$ integral \cite[Chapter 4]{DaPrato}. It is interpreted in the mild sense \cite[p. 161]{DaPrato}, meaning that an $H$-valued process $(U(t))_{t\in [0,T]}$ is a solution of \eqref{spde} if 
\begin{align}\label{mildsolution}
 U(t) = \e^{At}u_0 + \int_0^t \e^{A(t-s)} F(U(s))\; \de s + \int_0^t \e^{A(t-s)}\; \de W^Q(s).
\end{align}
We now give some details on the properties of each component of the SPDE \eqref{spde}.\\
\underline{The linear sectorial operator $A=\Delta$.}\\
It is known that for Dirichlet, Neumann and mixed boundary conditions, the Laplacian has a complete set of orthonormal smooth eigenfunctions $e_i$, $i\in\mathbb{N}$. Also, the eigenvalues $-\lambda_i$ are all negative and real having exactly one accumulation point at $-\infty$ (see e.g. \cite[Appendix W]{music}; we use a notation where $\lambda_i>0$, $i\in \mathbb{N}$). Hence, it is possible to give a spectral representation
\begin{align*}
\textstyle
 \Delta f(x) = \sum_{i=1}^\infty -\lambda_i\langle f,e_i \rangle_{H} e_i(x)
\end{align*}
for $f\in D(\Delta)$. We assume Dirichlet boundary conditions, so we call the $\lambda_i$ Dirichlet eigenvalues and the $e_i$ Dirichlet eigenfunctions.\\
\underline{The nonlinear operator $F$.}\\
In our setting, the operator $F$ on $H$ is a Nemytskii operator, i.e. there is an associated function $ u_F: \mathbb{R} \rightarrow \mathbb{R}$ such that the operator $F$ can be defined as $[F(g)](x) = u_F(g(x))$ for $g\in H$.
For $F$, we adopt the assumption from \cite[Assumption 2.4]{KloedenRoyal} and require that $F$ be two times continuously Fr\'{e}chet differentiable, and that its derivatives satisfy
\begin{align}\label{nonlinearbeginning}
 \|F'(u)-F'(v)\| &\leq \tilde L\|u-v \|, \\
 \| (-\Delta)^{-r}F'(u) (-\Delta)^r w \| &\leq \tilde L \| w \|
\end{align}
for all $u,v\in H$ and $w\in D((-\Delta)^r)$ for $r\in \{0,\ 1/2,\ 1\}$ and
\begin{align}\label{nonlinearend}
 \| \Delta^{-1} F''(u)(v,w) \| \leq \tilde L \| (-\Delta)^{-1/2} v \|\cdot \| (-\Delta)^{-1/2} w \|
\end{align}
for all $u,v,w\in H$ and $\tilde L>0$ a constant.
\\
\underline{The Wiener process $W^Q$.}\\
The process $W^Q=(W^Q(t))_{t\in[0,T]}$ is a $Q$-cylindrical Wiener process \cite[\S 4.1.2]{DaPrato}\cite[\S 3.2.35]{Lototsky}, i.e. a collection
\begin{align*}
 \{W_u^Q(t),\ u\in H,\ t\in[0,T]\}
\end{align*}
of zero-mean Gaussian random variables that satisfy 
\begin{align*}
 \mathbb{E}(W_u^Q(t)W_v^Q(s))=\langle Q u,v \rangle_H\cdot \min(t,s)
\end{align*}
for $s,t\in[0,T]$, $u,v\in H$. Then, for an orthonormal basis $\{e^Q_n\}_{n\in\mathbb{N}}$ of eigenfunctions of $Q$ with corresponding eigenvalues $\{q_k\}_{k\in\mathbb{N}}$, the process $W^Q$ has the series expansion in $H$ \cite[Proposition 4.3]{DaPrato}
\begin{align}\label{wienerseriesexpansion}
\textstyle
 W^Q(t) = \sum_{n=1}^\infty \sqrt{q_n} \beta_n(t) e_n^Q,
\end{align}
where $\beta_i$, $i\in\mathbb{N}$, are independent real-valued standard Brownian motions. We assume that $Q$ and the operator $\Delta$ both have a joint eigenbasis consisting of the same functions, i.e. we assume $\{e_n^Q\}_{n\in\mathbb{N}}=\{e_n\}_{n\in\mathbb{N}}$ (for details on this assumption, see Appendix \ref{appendixa}).  The stochastic integral with respect to the cylindrical $Q$-Wiener process is defined in a similar way to the classical It\^{o} integral (for details, see \cite[2.5]{LiuRoeckner}).\\
In the literature on SPDEs, at least one of the two constraints have to be made: Either a one-dimensional spatial domain is considered, where white noise as well as trace-class noise can be considered, or the covariance operator $Q$ cannot be allowed to be $Q=I$, i.e. space-time white noise cannot be considered. This is due to the fact that a solution of an SPDE on a multidimensional domain (i.e. $D\subset \R^n$, $n\geq 2$) with space-time white noise (i.e. $W^Q$ is a cylindrical Wiener process with $Q=I$) does not exist even for $F=0$, since vital regularity properties of the stochastic convolution
\begin{align*}
\textstyle
 W_A^Q(t) = \int_0^t \e^{A(t-s)}\; \de W^Q(s),\quad t\in [0,T],
\end{align*}
in \eqref{mildsolution}
are violated. See \cite[Example 5.7]{DaPrato}  \cite[p. 329]{Walsh} for discussions of these failures in spatial dimensions of two and higher. The deciding factor is the asymptotic behavior of Dirichlet eigenvalues, which has been described by Weyl's law \cite[4.1]{Grebenkov}: If $0<\lambda_1< \lambda_2\leq \lambda_3 \leq \dots$ are the Dirichlet eigenvalues of the negative Laplacian for an arbitrary domain $D\subset \R^d$, then
\begin{align}\label{weylslaw}
\textstyle
 \lambda_n \propto \frac{4\pi^2}{(\omega_d \mu_d(D))^{2/d}}n^{2/d},\quad n\rightarrow \infty,
\end{align}
where $\mu_d(D)$ is the $d$-dimensional Lebesgue measure of $D$ and 
$\omega_d$ is the volume of the $d$-dimensional unit ball. The regularity condition \cite[Assumption 2.2]{KloedenRoyal} for the exponential Euler scheme which we will use requires
\begin{align}\label{qncondition}
 \sum_{n=1}^\infty (\lambda_n)^{2\gamma-1} q_n <\infty
\end{align}
for some $\gamma\in(0,1)$, which for space-time white noise, i.e. $q_n=1$ for all $n\in\mathbb{N}$, is not fulfilled in two dimensions. Instead, we assume that $q_n = \mathcal{O}(n^{-\alpha})$ for some $\alpha>0$, in which case \eqref{qncondition} holds if $\gamma <\alpha/2$. In our numerical experiments, we will use $q_n = n^{-2+\epsilon_q}$ for an arbitrarily small $\epsilon_q>0$ so that $\gamma$ can be arbitrarily close to $1$ (adding $\epsilon_q$ is mostly a formal requirement which is justified in Appendix \ref{appendixa}).
\section{Numerical scheme.}\label{section_Numericalscheme}
In this work, we combine two different numerical methods: The spectral Galerkin exponential Euler scheme is used in two spatial dimensions with  eigenvalues and corresponding eigenfunctions which are numerically approximated using a boundary element collocation method in combination with the Beyn algorithm.\\ We use the spectral Galerkin method for approximation in space and the exponential Euler method for approximation in time to numerically solve the SPDE \eqref{spde} (Section \ref{exponentialEuler}). The spectral Galerkin scheme relies on projecting 
the values of the process from $H$ to the $N$-dimensional subspace $\text{span}(e_1,\dots,e_N)$, where the $e_i$ are the Dirichlet eigenfunctions of $\Delta$ (spatial approximation, Sections \ref{spatialapproximation}, \ref{computingonb}), and
then discretizing the mild solution \eqref{mildsolution} in time (temporal approximation, Section \ref{temporalapproximation}). 
If the SPDE has to be solved on a bounded domain $D\subset \R^2$, the Dirichlet eigenfunctions $e_i$ (and also the eigenvalues) for the domain $ D$ are needed. 
They will be computed using a boundary integral equation method (Section \ref{boundaryintegral}) with the Beyn integral algorithm (Section \ref{beyn}) as described in \cite{Kleefeldtransmission, Kleefeldshape}. Finally, we make a few comments on assembling an accurate orthonormal basis of Dirichlet eigenfunctions on a given shape (Section \ref{onb}).
\subsection{The spectral Galerkin exponential Euler scheme.}\label{exponentialEuler}
We recall the spectral Galerkin method and the exponential Euler scheme.
\subsubsection{Approximation in space.}\label{spatialapproximation}
We use a projection 
\begin{align*}
 P_N:H&\rightarrow H_N := \text{span}(e_1,\dots,e_N)\\
 f = \sum_{i=1}^\infty \langle f,e_i \rangle_{H} e_i &\mapsto \sum_{i=1}^N \langle f,e_i \rangle_{H} e_i
\end{align*}
with the inner product
\begin{align*}
 \langle f,e_i \rangle_{H} = \int_{ D} f(x) \cdot e_i(x)\;\de x
\end{align*}
and denote $U^N:=P_NU$, $F_N := P_NF$, $A_N:=P_NA$, $u_0^N:= P_N u_0$ and 
\begin{align*}
 W^Q_N := P_NW^Q = \sum_{i=1}^N \sqrt{q_i} e_i \beta_i(t).
\end{align*}
The projected SPDE
\begin{align*}
 \de U^N(t) = (A_N U^N(t) + F_N(U^N(t)))\de t + \de W_N^Q(t)
\end{align*}
has the mild solution satisfying for $t\in [0,T]$
\begin{align}\label{mildsolution2}
 U^N(t) = \e^{A_N t}u_0^N + \int_0^t \e^{A_N(t-s)}F_N(U^N(s))\;\de s + \int_0^t \e^{A_N(t-s)}\;\de W_N^Q(s).
\end{align}
\subsubsection{Approximation in time.}\label{temporalapproximation}
We introduce the equidistant time steps $t_i=h\cdot i$, $i=0,\dots,M$, with $h=T/M$, and denote $u_k^{N,M}=U^N(t_k)$. Equation \eqref{mildsolution2} can be rewritten in one time step from $t_k$ to $t_{k+1}$ as 
{\small
\begin{align*}%\label{timeapproximation1}
 u_{k+1}^{N,M} = \e^{A_N h} u_k^{N,M} + \int_{t_k}^{t_{k+1}} \e^{A_N(t_{k+1}-s)}F_N(u^{N,M}(s))\de s + \int_{t_k}^{t_{k+1}} \e^{A_N(t_{k+1}-s)}\de W_N^Q(s).
\end{align*}
}%
We now make the approximation $F_N(u^{N,M}(s)) \approx F_N(u^{N,M}(t_k))$.
Let $\widetilde{U^N}$ be the process after this approximation and let $v_k^{N,M}:=\widetilde{U^N}(t_k)$. We denote the inner products by
\begin{align*}%\label{innerproducts}
v_{k,j}^{N,M}:=\langle v_k^{N,M},e_j \rangle_{H},\quad F_N^j(v_k^{N,M}):=\langle F_N(v_k^{N,M}),e_j \rangle_{H}.
\end{align*}
Taking the inner product with $e_j$ on both sides and applying the above approximation, we obtain for the approximate Fourier-Galerkin coefficients
\begin{align*}
 v^{N,M}_{k+1,j} 
 & = \langle \e^{A_N h} v_k^{N,M},e_j \rangle + \langle \int_{t_k}^{t_{k+1}}\e^{A_N(t_{k+1}-s)}F_N(v_k^{N,M})\; \de s,e_j \rangle \\
& \quad + \langle \int_{t_k}^{t_{k+1}} \e^{A_N(t_{k+1}-s)}\; \de W_N^Q(s), e_j \rangle \\
 & = \e^{-\lambda_j h} v_{k,j}^{N,M} + F_N^j(v_k^{N,M})\cdot\int_{t_k}^{t_{k+1}} \e^{-\lambda_j(t_{k+1}-s)}\; \de s \\
 & \quad + \int_{t_k}^{t_{k+1}} \e^{-\lambda_j (t_{k+1}-s)} \sqrt{q_j}\; \de \beta_j(s),
\end{align*}
using linearity of the inner product and the fact that $ A_N h e_i = -\lambda_i h e_i$ and hence $\e^{A_N h} e_i = \e^{-\lambda_i h} e_i$ for $h>0$. The $\de \beta_j(s)$ integral is normally distributed with mean zero and its variance can be computed using the It\^{o} isometry. We finally obtain
\begin{align}\label{finalscheme}
 v_{k+1,j}^{N,M} = \e^{-\lambda_j h} v_{k,j}^{N,M} + \frac{1-\e^{-\lambda_j h}}{\lambda_j}F_N^j(v_k^{N,M}) + \left(\frac{q_j}{2\lambda_j}(1-\e^{-2\lambda_j h})\right)^{1/2} R_k^j,
\end{align}
where the $R_k^j$ are independent identically distributed standard normally distributed random variables, $j=1,\dots,N$, $k=0,\dots,M-1$. This is the exponential Euler scheme, as described in \cite[Section 3]{KloedenRoyal}.
\subsection{The boundary integral equation method.}\label{boundaryintegral}
The boundary integral equation method aims to solve boundary value problems such as the Helmholtz equation with Dirichlet or Neumann boundary conditions by writing the solution as an integral layer operator involving an unknown density, approximating this density and subsequently reinserting it into the integral to compute the solution.\\
Let $\Phi_\kappa(x,y)$ be the fundamental solution 
\begin{align}
\Phi_\kappa(x,y) = \text{i} H_0^{(1)}(\kappa \| x-y\|)/4,\quad x\neq y
\end{align}
of the Helmholtz equation
\begin{align}\label{helmholtz}
 -\Delta u + \kappa^2 u &= 0\;\; \text{in}\  D,\qquad  u = 0\;\; \text{on}\ \partial D,
\end{align}
where $H_0^{(1)}$ denotes the first-kind Hankel function of order zero (see \cite[p. 66]{colton2013inverse}) and let $\kappa\in\mathbb{C}\backslash\{ 0 \}$ be fixed.
For a point $y$ on the boundary $\partial D$, we denote the exterior normal by $\nu(y)$ and we use the double layer potential
\begin{align}\label{dlpotential}
 \DL_\kappa[\psi](x) = \int_\Gamma \psi(y)\cdot \partial_{\nu(y)}\Phi_\kappa (x,y)\;\de s(y),\quad x\in \R^2\backslash\partial D,
\end{align}
where $\psi\in C(\partial D)$ is the unknown density function. We also define
\begin{align*}
 \text{D}_\kappa[\psi](x) = \int_\Gamma \psi(y)\cdot \partial_{\nu(y)}\Phi_\kappa(x,y)\;\de s(y),\quad x\in\partial D.
\end{align*}
For $x\in \partial D$, it is known \cite[Theorem 3.1]{colton2013inverse} that the so-called jump relations
\begin{align}\begin{split}\label{jumprelation}
    \lim_{h\rightarrow 0} &\int_{\Gamma} \psi(y) \partial_{\nu(y)} \Phi_\kappa (x \pm h\nu(x),y)\; \de s(y) = \\ &\int_{\Gamma} \psi(y) \partial_{\nu(y)} \Phi_\kappa (x,y)\; \de s(y) \pm \frac{1}{2}\psi(x)
    \end{split}
\end{align}
hold. 
Since the fundamental solution $\Phi_\kappa$ satisfies the Helmholtz equation \eqref{helmholtz}, so does the double-layer potential \eqref{dlpotential}. Conversely, every solution of \eqref{helmholtz} can be written as a double-layer potential \cite[p. 39]{colton2013inverse}. We use the double layer ansatz
\begin{align}\label{dlansatz}
 u(x) = \DL_\kappa[\psi](x),\quad x\in D
\end{align}
which, after letting $x$ tend towards the boundary $\partial  D$, with the jump relation \eqref{jumprelation} becomes
\begin{align}\label{evproblem1}
 -\frac{1}{2} \psi(x) + \text{D}_\kappa[\psi](x) = 0,\quad x\in \partial D.
\end{align}
Writing $\mathrm{M}(\kappa):=-\frac{1}{2}\mathrm{I}+\mathrm{D}_\kappa$, \eqref{evproblem1} becomes $\mathrm{M}(\kappa)\psi = 0$, which is a nonlinear eigenvalue problem. First, the boundary is subdivided into $n_f$ curved elements called $\Delta_1,\dots,\Delta_{n_f}$, and \eqref{evproblem1} becomes
\begin{align}\label{discretization1}
 -\frac{1}{2}\psi(x) + \sum_{j=1}^{n_f} \int_{\Delta_j} \psi(y) \partial_{\nu(y)}\Phi_\kappa(x,y)\;\de s(y) = 0,\quad x\in \partial D.
\end{align}
Let $m_j:[0,1]\rightarrow \Delta_j$, $j=1,\dots, n_f$, be bijective, continuous maps. With a change of variables, \eqref{discretization1} becomes
\begin{align}\label{discretization2}
 -\frac{1}{2}\psi(x) + \sum_{j=1}^{n_f} \int_0^1 \psi(m_j(s))\partial_{\nu(m_j(s))}\Phi_\kappa(x,m_j(s)) \| \partial_s m_j(s)\| \;\de s(s) = 0,\quad x\in\partial D.
\end{align}
We now denote by $L_i(s)$ the quadratic Lagrange polynomials and by $\widehat{L_i}(s)$ the generalized Lagrange basis functions, $i=1,2,3$:
\begin{alignat*}{3}
 &L_1(s) = (1-s)(1-2s),\qquad  &&L_2(s) = 4s(1-s),\qquad\qquad\;\;\;   &&L_3(s) = s(2s-1), \\
 &\widehat{L_1}(s) = \frac{1-s-\alpha}{1-2\alpha}\frac{1-2s}{1-2\alpha}, &&\widehat{L_2}(s) = 4\frac{s-\alpha}{1-2\alpha}\frac{1-s-\alpha}{1-2\alpha}, &&\widehat{L_3}(s) = \frac{s\alpha}{1-2\alpha}\frac{2s-1}{1-2\alpha}
 \end{alignat*}
and we fix a parameter
\begin{align}
    \label{alpha_firstappearance}
    \alpha\in(0,1/2).
\end{align}
 %In \cite[Section 4.6]{hotspots}, it is explained that while different vales of $\alpha$ lead to optimal convergence rates for smooth boundaries, the choice $\alpha=(1-\sqrt{3/5})/2$ assured at least cubic convergence rate for a given eigenvalue in all computations, so we also choose $\alpha$ as above.
Each $m_j(s)$ is approximated by a quadratic interpolation polynomial 
\begin{align*}
 m_j(s)\approx \tilde m_j(s) = \sum_{i=1}^3 v^{(j)}_i L_i(s),
\end{align*}
where $v_i^{(j)}$, $i=1,3$ are the end points of the $j$-th boundary element and $v_2^{(j)}$ is its midpoint. We define the collocation nodes $\tilde v_{j,k}:=\tilde m_j(q_k)$ for $j=1,\dots,n_f$, $k=1,2,3$ and $q_1=\alpha,\ q_2=1/2,\ q_3 = 1-\alpha$. The unknown function $\psi(\tilde m_j(s))$ is approximated on each piece $\Delta_j$ by a quadratic interpolation polynomial 

    \begin{align*}
\sum_{k=1}^3 \psi(\tilde m_j(q_k))\widehat{ L_k}(s)=\sum_{k=1}^3 \psi(\tilde v_{j,k})\widehat{L_k}(s).        
    \end{align*}

 Equation \eqref{discretization2} then becomes
\begin{align*}
 -\frac{1}{2}\psi(x) + \sum_{j=1}^{n_f}\sum_{k=1}^3 \int_0^1 \partial_{\nu(\tilde m_j(s))} \Phi_\kappa(x,\tilde m_j(s))\| \partial_s \tilde m_j(s)\| \widehat{L_k}(s)\;\de s(s) \psi(\tilde v_{j,k}) = r(x),
\end{align*}
with a residue $r(x)$ which results from the various approximations. For the boundary element collocation method, we require $r(\tilde v_{i,\ell})=0$ \cite[p. 50]{atkinson} and obtain the linear system of size $3n_f\times 3n_f$
\begin{align}\label{discretization3}
 -\frac{1}{2}\psi(\tilde v_{i,\ell}) + \sum_{j=1}^{n_f}\sum_{k=1}^3 a_{i,j,k,\ell} \psi(\tilde v_{j,k}) = 0
\end{align}
with $i=1,\dots,n_f$, $\ell=1,2,3$ and with the integrals
\begin{align*}
 a_{i,j,k,\ell} = \int_0^1 \partial_{\nu(\tilde m_j(s))} \Phi_\kappa(\tilde v_{i,\ell},\tilde m_j(s))\| \partial_s \tilde m_j(s)\| \widehat{L_k}(s)\;\de s(s)
\end{align*}
which are approximated using a Gau{\ss}-Kronrod quadrature. We write \eqref{discretization3} abstractly as the nonlinear eigenvalue problem $M(\kappa)\phi = 0$,
where
\begin{align}\label{approximatedensity}
 \phi = (\psi(v_{1,1}),\psi(v_{1,2}),\psi(v_{1,3}),\psi(v_{2,1}),\dots,\psi(v_{n_f,3}))^\top.
\end{align}
\subsection{The Beyn integral algorithm.}\label{beyn}
For a given wavenumber $\kappa$, the solution of \eqref{helmholtz} can be found by inserting the obtained density $\psi$ back into \eqref{dlansatz}, so the eigenfunction value can be 
computed at any given point inside $ D$. But first, the nonlinear eigenvalue problem $M(\kappa)\phi = 0$ for the density must be solved. This is accomplished with the Beyn algorithm 
\cite[p. 3849]{beyn}, which makes use of Keldysh's theorem and is based on complex contour integrals of the resolvent $M^{-1}(\kappa)$, whose poles are the Dirichlet eigenvalues. 
\begin{theorem}[Keldysh's Theorem.]
 Let $D\subset \mathbb{C}$ be a domain, and let $v^* = \overline{v^\top}$ denote the conjugate transpose of a vector or matrix. Suppose that $C\subset D$ is a compact subset containing only simple eigenvalues $\lambda_n$, $n=1,\dots,n(\gamma)$ with eigenvectors $v_n$, $w_n$ satisfying 
 \begin{align*}
  M(\lambda_n)v_n = 0,\quad w_n^* M(\lambda_n) = 0,\quad w_n^*M'(\lambda_n) v_n = 1.
 \end{align*}
Then there exists some neighborhood $U$ of $C$ in $D$ and a holomorphic function $R:U\rightarrow \mathbb{C}^{3n_f\times 3n_f}$ such that
\begin{align}\label{keldysh}
 M(z)^{-1} = \sum_{n=1}^{n(\gamma)} \frac{1}{z-\lambda_n}v_nw_n^* + R(z),\quad z\in U\backslash \{\lambda_1,\dots,\lambda_{n(\gamma)}\}.
\end{align}
Let $\sigma(M)$ denote the spectrum of $M$ and let $\gamma\subset D$ be a contour which satisfies $\gamma\cap \sigma(M) = \emptyset$ and containing $n(\gamma)$ eigenvalues $\lambda_1,\dots,\lambda_{n(\gamma)}$.
Let $f:D\rightarrow \mathbb{C}$ be holomorphic. Then, applying the residue theorem to \eqref{keldysh} we have
\begin{align*}
 \frac{1}{2\pi \mathrm{i}}\int_\gamma f(z) M(z)^{-1}\;\de z = \sum_{n=1}^{n(\gamma)} f(\lambda_n) v_n w_n^*.
\end{align*}
\end{theorem}
For a specified contour $\gamma(t) = \mu+R \cos(t)+ \text{i} R\sin(t)$, $\mu\in\R$, $t\in[0,2\pi)$, $R\in\mathbb{R}_{>0}$ and the number $n(\gamma)$ of eigenvalues enclosed in $\gamma$ (with multiplicity), we pick
a random matrix $\hat{V}\in \mathbb{C}^{3n_f\times \ell}$ with $3n_f\gg \ell\geq n(\gamma)$. We compute the contour integrals
\begin{alignat*}{2}
\textstyle
 A_0 &= \frac{1}{2\pi \text{i}}\int_\gamma M^{-1}(\kappa)\hat {V}\;\de s(\kappa) &&= \frac{1}{2\pi\text{i}}\int_0^{2\pi} M^{-1}(\gamma(t))\hat{{V}} \gamma'(t)\;\de s(t) \\
 &= \sum_{n=1}^{n(\gamma)} v_n w_n^*\hat V = VW^*\hat V, \\
 A_1 &= \frac{1}{2\pi \text{i}}\int_\gamma\kappa M^{-1}(\kappa)\hat {V}\;\de s(\kappa) &&= \frac{1}{2\pi\text{i}}\int_0^{2\pi}\gamma(t) M^{-1}(\gamma(t))\hat{{V}} \gamma'(t)\;\de s(t) \\
 &= \sum_{n=1}^{n(\gamma)} \lambda_n v_n w_n^* \hat V = V\Lambda W^*\hat V,
\end{alignat*}
with $\Lambda = \text{diag}(\lambda_1,\dots,\lambda_{n(\gamma)})$. These integrals are approximated using the trapezoidal rule:
\begin{align}\label{trapezoidalrule}
 A_{0,K} &= \frac{1}{\text{i}K}\sum_{j=0}^{K-1} M^{-1}(\gamma(t_j))\hat{V} \gamma'(t_j),\; 
 A_{1,K} = \frac{1}{\text{i}K}\sum_{j=0}^{K-1} \gamma(t_j) M^{-1}(\gamma(t_j))\hat{V} \gamma'(t_j),
\end{align}
with $t_j=\frac{2\pi j}{K},\ j=0,\dots,K$. Next, we compute a singular value decomposition $A_{0,K} = V\Sigma W^*$ with $V\in \mathbb{C}^{3n_f\times \ell}$, $\Sigma\in \mathbb{C}^{\ell\times \ell}$ and $W\in \mathbb{C}^{\ell\times \ell}$. For $\Sigma=\text{diag}(\sigma_1,\sigma_2,\dots,\sigma_\ell)$ and a given tolerance $\epsilon_{\mathrm{sing}}$, we find $n(\gamma)$ such that
\begin{align}\label{esingline}
 \sigma_1\geq \dots \geq \sigma_{n(\gamma)}>\epsilon_{\mathrm{sing}}>\sigma_{n(\gamma)+1}\geq \dots \geq \sigma_\ell.
\end{align}
The singular values $\sigma_{n(\gamma)+1}\geq \dots \geq \sigma_\ell$ which are below the tolerance are likely to be different from zero only due to rounding errors and are discarded. With the matrices $V_0:= (V_{ij})_{1\leq i\leq m,1\leq j\leq n(\gamma)},\ \Sigma_0 := (\Sigma_{ij})_{1\leq i\leq n(\gamma),1\leq j\leq n(\gamma)}$ as well as $W_0 := (W_{ij})_{1\leq i\leq \ell,1\leq j\leq n(\gamma)}$
there is a regular matrix $S\in \mathbb{C}^{n(\gamma)\times n(\gamma)}$ such that $V=V_0 S$, and after rewriting some matrices, we have that
\begin{align*}
 S\Lambda S^{-1} = V_0^* A_1 W_0 \Sigma_0^{-1}.
\end{align*}
The matrix on the right-hand side therefore has the same eigenvalues and eigenvectors as $\Lambda$. We compute the $n(\gamma)$ eigenvalues and 
eigenvectors $s_i$ of the matrix $V_0^* A_{1,N}W_0 \Sigma_0^{-1}\in \mathbb{C}^{n(\gamma)\times n(\gamma)}$. The $i$-th eigenfunction $u_i$ is then approximated by inserting $\phi = V_0 s_i$ (see \eqref{approximatedensity}) into
\begin{align*}
 u_i(x) = \mathrm{DL}_\kappa [\psi](x) \approx \sum_{j=1}^{n_f}\sum_{k=1}^3 \hat a_{j,k}(x)\cdot \psi(\tilde v_{j,k}),
\end{align*}
where
\begin{align*}
 \hat a_{j,k}(x) = \int_0^1 \partial_{\nu(\tilde m_j(s))}\Phi_\kappa(x,\tilde m_j(s))\|\partial_s \tilde m_j(s)\| \widehat{L_k}(s)\; \de s(s),
\end{align*}
which again are approximated using a Gau{\ss}-Kronrod quadrature.

\begin{figure}
\newcommand{\Height}{5 cm}
\newcommand{\Width}{5 cm}
\newcommand{\X}{2.5}
\newcommand{\Yone}{1.0}
\newcommand{\Ytwo}{0.5}
\centering
 \begin{tabular}{@{}c@{}}
 % 1
\begin{tikzpicture}
	\begin{axis}[
		height=\Height, 
		width=\Width,
		axis x line=bottom,
		axis y line=left,
		xlabel = {$\kappa$},
		ylabel= {$n_f$},
		clip=false,
		ymin=100,ymax=520,
        xmin=0,xmax=85,
		xtick={0,10,20,30,40,50,60,70,80},
		ytick={100,150,200,250,300,350,400,450,500},
		tick label style={font=\footnotesize},
		tick align=outside
		]
		\addplot[color=green,only marks,mark=o] table {data11g.txt};
		%\addplot[color=red,only marks,mark=o] table {data11r.txt};
	\end{axis}
	\node[color=black] at (\X, \Yone) (p1) {\footnotesize{$\epsilon^{(\lambda)}=2\cdotp 10^{-4}$}};
	\node[color=black] at (\X, \Ytwo) (p1) {\footnotesize{$\epsilon^{(\eta)}=5\cdotp 10^{-6}$}};
\end{tikzpicture}
% 2
\begin{tikzpicture}
	\begin{axis}[
		height=\Height, 
		width=\Width,
		axis x line=bottom,
		axis y line=left,
		xlabel = {$\kappa$},
		ylabel= {$n_f$},
		clip=false,
		ymin=100,ymax=520,
        xmin=0,xmax=85,
		xtick={0,10,20,30,40,50,60,70,80},
		ytick={100,150,200,250,300,350,400,450,500},
		tick label style={font=\footnotesize},
		tick align=outside
		]
		\addplot[color=green,only marks,mark=o] table {data21g.txt};
		%\addplot[color=red,only marks,mark=o] table {data21r.txt};
		\addplot[dashed, color=blue] table {bemslope1.txt};	
	\end{axis}
	\node[color=black] at (\X, \Yone) (p1) {\footnotesize{$\epsilon^{(\lambda)}=2\cdotp 10^{-4}$}};
	\node[color=black] at (\X, \Ytwo) (p1) {\footnotesize{$\epsilon^{(\eta)}=5\cdotp 10^{-6}$}};
\end{tikzpicture}
\\
% 3
\begin{tikzpicture}
	\begin{axis}[
		height=\Height, 
		width=\Width,
		axis x line=bottom,
		axis y line=left,
		xlabel = {$\kappa$},
		ylabel= {$n_f$},
		clip=false,
        ymin=100,ymax=520,
        xmin=0,xmax=85,
		xtick={0,10,20,30,40,50,60,70,80},
		ytick={100,150,200,250,300,350,400,450,500},
		tick label style={font=\footnotesize},
		tick align=outside
		]
		\addplot[color=green,only marks,mark=o] table {data12g.txt};
		\addplot[color=red,only marks,mark=o] table {data12r.txt};
	\end{axis}
	\node[color=black] at (\X, \Yone) (p1) {\footnotesize{$\epsilon^{(\lambda)}=1\cdotp 10^{-4}$}};
	\node[ color=black] at (\X, \Ytwo) (p1) {\footnotesize{$\epsilon^{(\eta)}=5\cdotp 10^{-6}$}};
\end{tikzpicture}
% 4
\begin{tikzpicture}
	\begin{axis}[
		height=\Height, 
		width=\Width,
		axis x line=bottom,
		axis y line=left,
		xlabel = {$\kappa$},
		ylabel= {$n_f$},
		clip=false,
		ymin=100,ymax=520,
        xmin=0,xmax=85,
		xtick={0,10,20,30,40,50,60,70,80},
		ytick={100,150,200,250,300,350,400,450,500},
		tick label style={font=\footnotesize},
		tick align=outside
		]
		\addplot[color=green,only marks,mark=o] table {data22g.txt};
		\addplot[color=red,only marks,mark=o] table {data22r.txt};
		\addplot[dashed, color=blue] table {bemslope2.txt};	
	\end{axis}
	\node[color=black] at (\X, \Yone) (p1) {\footnotesize{$\epsilon^{(\lambda)}=1\cdotp 10^{-4}$}};
	\node[color=black] at (\X, \Ytwo) (p1) {\footnotesize{$\epsilon^{(\eta)}=5\cdotp 10^{-6}$}};
\end{tikzpicture}
\\
% 5
\begin{tikzpicture}
	\begin{axis}[
		height=\Height, 
		width=\Width,
		axis x line=bottom,
		axis y line=left,
		xlabel = {$\kappa$},
		ylabel= {$n_f$},
		clip=false,
		ymin=100,ymax=520,
        xmin=0,xmax=85,
		xtick={0,10,20,30,40,50,60,70,80},
		ytick={100,150,200,250,300,350,400,450,500},
		tick label style={font=\footnotesize},
		tick align=outside
		]
		\addplot[color=green,only marks,mark=o] table {data13g.txt};
		\addplot[color=red,only marks,mark=o] table {data13r.txt};
	\end{axis}
	\node[color=black] at (\X, \Yone) (p1) {\footnotesize{$\epsilon^{(\lambda)}=2\cdotp 10^{-4}$}};
	\node[color=black] at (\X, \Ytwo) (p1) {\footnotesize{$\epsilon^{(\eta)}=2\cdotp 10^{-6}$}};
\end{tikzpicture}
% 6
\begin{tikzpicture}
	\begin{axis}[
		height=\Height, 
		width=\Width,
		axis x line=bottom,
		axis y line=left,
		xlabel = {$\kappa$},
		ylabel= {$n_f$},
		clip=false,
		ymin=100,ymax=520,
		xmin=0,xmax=85,
		xtick={0,10,20,30,40,50,60,70,80},
		ytick={100,150,200,250,300,350,400,450,500},
		tick label style={font=\footnotesize},
		tick align=outside
		]
		\addplot[color=green,only marks,mark=o] table {data23g.txt};
		\addplot[color=red,only marks,mark=o] table {data23r.txt};
		\addplot[dashed, color=blue] table {bemslope3.txt};		
	\end{axis}
	\node[color=black] at (\X, \Yone) (p1) {\footnotesize{$\epsilon^{(\lambda)}=2\cdotp 10^{-4}$}};
	\node[color=black] at (\X, \Ytwo) (p1) {\footnotesize{$\epsilon^{(\eta)}=2\cdotp 10^{-6}$}};
\end{tikzpicture}
\end{tabular}
 \caption{The number of boundary elements to approximate wavenumbers and their respective eigenfunctions to stay below certain error tolerances $\epsilon^{(\lambda)}$ and
 $\epsilon^{(\eta)}$. For each base function, we computed a reference solution with $n_f^{\mathrm{ref}} = 600$ boundary elements and compared the same function with less boundary elements to it, with $n_f\in\{50,75,100,150,200,300,400,500\}$ boundary elements used. The plots on the left show raw data, with the least number of required boundary elements. The plots on the right
 show linearly interpolated numbers of needed boundary elements based off the data from the plots on the left. We can observe the linear relationship between wavenumber and needed boundary elements $n_f$, and that the needed $n_f$ increases as the error tolerances decrease (in red are the values which exceed the error tolerance for every tested $n_f$). The blue reference slopes are functions $f(x)=mx+b$, for $m\in\{7, 8.5, 10 \}$, $b\in \{80, 65, 70\}$, from top to bottom.}
 \label{boundaryelements}
\end{figure}
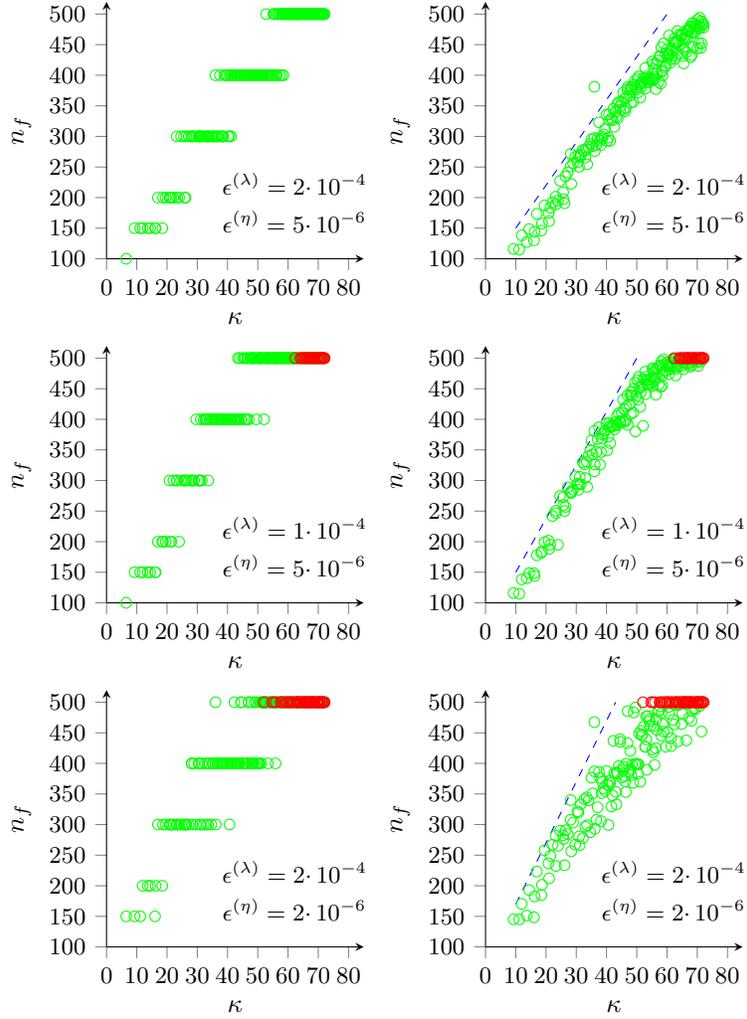

\subsection{Preprocessing: Building an orthonormal basis (ONB)}\label{onb}%\hfill \\
The analytical expressions of the Dirichlet eigenfunctions are known for some shapes. In the following, we will consider a shape which we call the Peanut shape (see Figure \ref{Reference1200600}), whose Dirichlet eigenvalues and eigenfunctions are not known analytically and whose parametrization for $t\in [0,2\pi)$ is given by 
{\small
\begin{equation*}
\begin{split}
\textstyle
 \left(\begin{array}{c} 0.06\cdot((\cos(t)+2)\cdot(\cos(t+0.6)+2)\cdot(0.1\cdot\cos(3t)+2))-0.1 \\
                      0.06\cdot(\sin(t)+2)\cdot (\sin(t-0.5)+2)\cdot (0.4\cdot \cos(2t)+2)\cdot(0.1\cdot\sin(4t)+1))-0.06
                      \end{array}\right).
\end{split}
\end{equation*}
}%
For the boundary element collocation method, we choose $\alpha = (1-\sqrt{3/5})/2$ (see \eqref{alpha_firstappearance}). Note that while any $\alpha \in (0,\frac{1}{2})$ ensures at least cubic convergence, empirically it was seen \cite[Section 4.6]{hotspots} that the choice $\alpha = (1-\sqrt{3/5})/2$ yielded superconvergence for domains with a $\mathcal{C}^2$ boundary.
Whenever the Beyn integral algorithm is used to compute its Dirichlet eigenvalues and eigenfunctions, we use $\epsilon_{\mathrm{sing}}=10^{-4}$ as singular value tolerance (see \eqref{esingline}) and $K=24$ for the periodic trapezoidal rule (see \eqref{trapezoidalrule}).
\subsubsection{Convergence and accuracy of base functions.}\label{eigenfunctionconvergence} 
In previous work, the Beyn integral algorithm was used to compute Neumann eigenvalues, and a convergence with rates up to three for a fixed wavenumber could be observed \cite{hotspots}.
However, the eigenvalues and eigenfunctions become less accurate for higher wavenumbers, so that more boundary elements are needed for higher wavenumbers. Previous publications from high-frequency 
scattering theory using boundary element collocation methods to solve the Helmholtz equation have pointed out that this effort rises linearly with the wavenumber \cite{highfrequencyscattering,chandlerhighfrequency},  which we could confirm for the case of the Dirichlet eigenfunctions and eigenvalues (see Figure \ref{boundaryelements}). For the Peanut shape, we have conducted a convergence analysis for the first 200 eigenfunctions. For each wavenumber, we computed a reference function for $600$ boundary elements which we compare to functions with lower accuracy (computed with 500, 400, 300, 200, 150, 100, 75 and 50 boundary elements, 
respectively). For speeding up the computations, these functions were stored only on an equidistant $81\times 81$ grid inside $[0,1]^2$, while the final functions making up the ONB are stored on a $301\times 301$ grid to obtain results with better resolution (for the error analysis in Section \ref{section_erroranalysis}, we will denote the grid resolution as $R$ and will use $R=301$). The pointwise storage on the equidistant grid means that the inner products, i.e. the integrals in the scheme \eqref{finalscheme} can be approximated numerically. As it is easy to implement and has a convergence order of four \cite[p. 21 f.]{numint} (and hence, a convergence order of two for two-dimensional functions), we use the composite two-dimensional Simpson rule and set the numerically integrated function to be zero outside the boundary, so that the Simpson rule can be applied on the unit square. The application of the Simpson rule is possible since eigenfunctions of the Laplacian are always smooth \cite[\S 2.(i)]{Grebenkov}.\\
Suppose that we use $n_f^{(i)}$ boundary elements for eigenfunction number $i$, $i=1,\dots,N$. We fix error tolerances
\begin{align*}
  \epsilon^{(\lambda)} := \max_{i=1,\dots,N} |\widetilde{ \lambda_i}^{(n_f^{(i)})} - \tilde \lambda_i^{(600)}|,\ \epsilon^{(\eta)} &:= \max_{i=1,\dots,N} \|\tilde e_i^{(n_f^{(i)})} - \tilde e_i^{(600)}\|_{H},
\end{align*}
where the tilde denotes approximation and the superscript denotes that the approximation is done with $n_f^{(i)}$ boundary elements. Depending on which $\epsilon^{(\lambda)}$ and $\epsilon^{(\eta)}$ we choose, different wavenumber dependent amounts of boundary elements are needed (see Figure \ref{boundaryelements}). We will later observe (see Example \ref{errorboundexample}) that for the error bound derived in Theorem \ref{strongerrortheorem} it is crucial to reduce $\epsilon^{(\lambda)}$ and $\epsilon^{(\eta)}$ as much as possible.
\subsubsection{Computation of the orthonormal basis.}\label{computingonb}
\begin{figure}
\definecolor{light_grey}{RGB}{140,140,140}
\newcommand\lambdaposa{0.8}
\newcommand\lambdaposb{2.5}
\newcommand\lambdaposc{4.55}
\newcommand\shift{7}
\newcommand\circleaa{0.8}
\newcommand\circleab{2.1}
\newcommand\circleac{3.5}
\newcommand\circlead{4.9}
\newcommand\circleradi{0.8}
    \centering
    \begin{tikzpicture}[scale=0.92]
    \draw [|-stealth] (0,0) -- (6,0);
    \draw [-stealth] (0,-1) -- (0,1);
    \node at (6.2,0) {$\mathbb{R}$};
    \node at (0.5,1) {$\mathbb{C}$};
    \node at (-0.25,0) {$0$};
    \node at (\lambdaposa,0.3) {\textcolor{light_grey}{$\lambda_1$}};
    \node at (\lambdaposb,0.3) {\textcolor{light_grey}{$\lambda_2$}};
    \node at (\lambdaposc,0.3) {\textcolor{light_grey}{$\lambda_3$}};
    \node at (5.1,0.3) {$\dots$};
    \node at (\lambdaposa,0) {\textcolor{light_grey}{$\times$}};
    \node at (\lambdaposb,0) {\textcolor{light_grey}{$\times$}};
    \node at (\lambdaposc,0) {\textcolor{light_grey}{$\times$}};
    \foreach \x in {\circleaa,\circleab,\circleac,\circlead}
    {
    \draw[->] (\x,0.8) arc [
    start angle=92, end angle=448,
    radius=\circleradi cm];
    }
    \node at (\circleaa,-1) {$\hat\gamma_1$};
    \node at (\circleab,-1) {$\hat\gamma_2$};
    \node at (\circleac,-1) {$\hat\gamma_3$};
    \node at (\circlead,-1) {$\hat\gamma_4$};
    \draw [|-stealth] (\shift,0) -- (6+\shift,0);
    \draw [-stealth] (\shift,-1) -- (\shift,1);
    \node at (\shift+6.2,0) {$\mathbb{R}$};
    \node at (\shift+0.5,1) {$\mathbb{C}$};
    \node at (\shift-0.25,0) {$0$};
    \node at (\shift+\lambdaposa,0.3) {{$\lambda_1$}};
    \node at (\shift+\lambdaposb,0.3) {{$\lambda_2$}};
    \node at (\shift+\lambdaposc,0.3) {{$\lambda_3$}};
    \node at (\shift+5.1,0.3) {$\dots$};
    \node at (\shift+\lambdaposa,0) {{$\times$}};
    \node at (\shift+\lambdaposb,0) {{$\times$}};
    \node at (\shift+\lambdaposc,0) {{$\times$}};
    \draw[->] (\shift+\lambdaposa,0.6) arc [
    start angle=92, end angle=448,
    radius=0.6cm];
    \draw[->] (\shift+2.5,0.7) arc [
    start angle=92, end angle=448,
    radius=0.7cm];
    \draw[->] (\shift+4.55,0.8) arc [
    start angle=92, end angle=448,
    radius=0.8cm];
    \node at (\shift+\lambdaposa,-0.8) {$\gamma_1$};
    \node at (\shift+2.5,-0.9) {$\gamma_2$};
    \node at (\shift+4.55,-1) {$\gamma_3$};
    \end{tikzpicture}
    \caption{In the first step, the eigenvalues are still unknown, so the real axis up to a certain point is covered with contours $\hat \gamma_i$ which cover the real axis up to the point where the desired number $N$ of eigenvalues has been found. For the subsequent, more precise computation, contours $\gamma_i$, $i=1,\dots, N$ around the found eigenvalues are chosen.}
    \label{contours}
\end{figure}
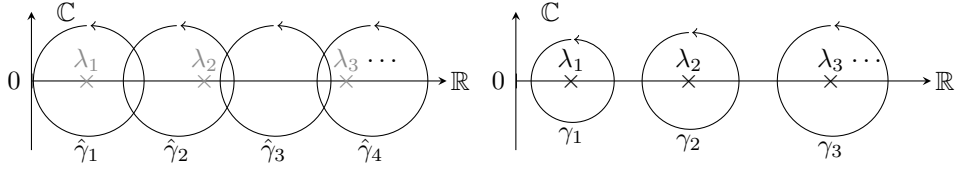

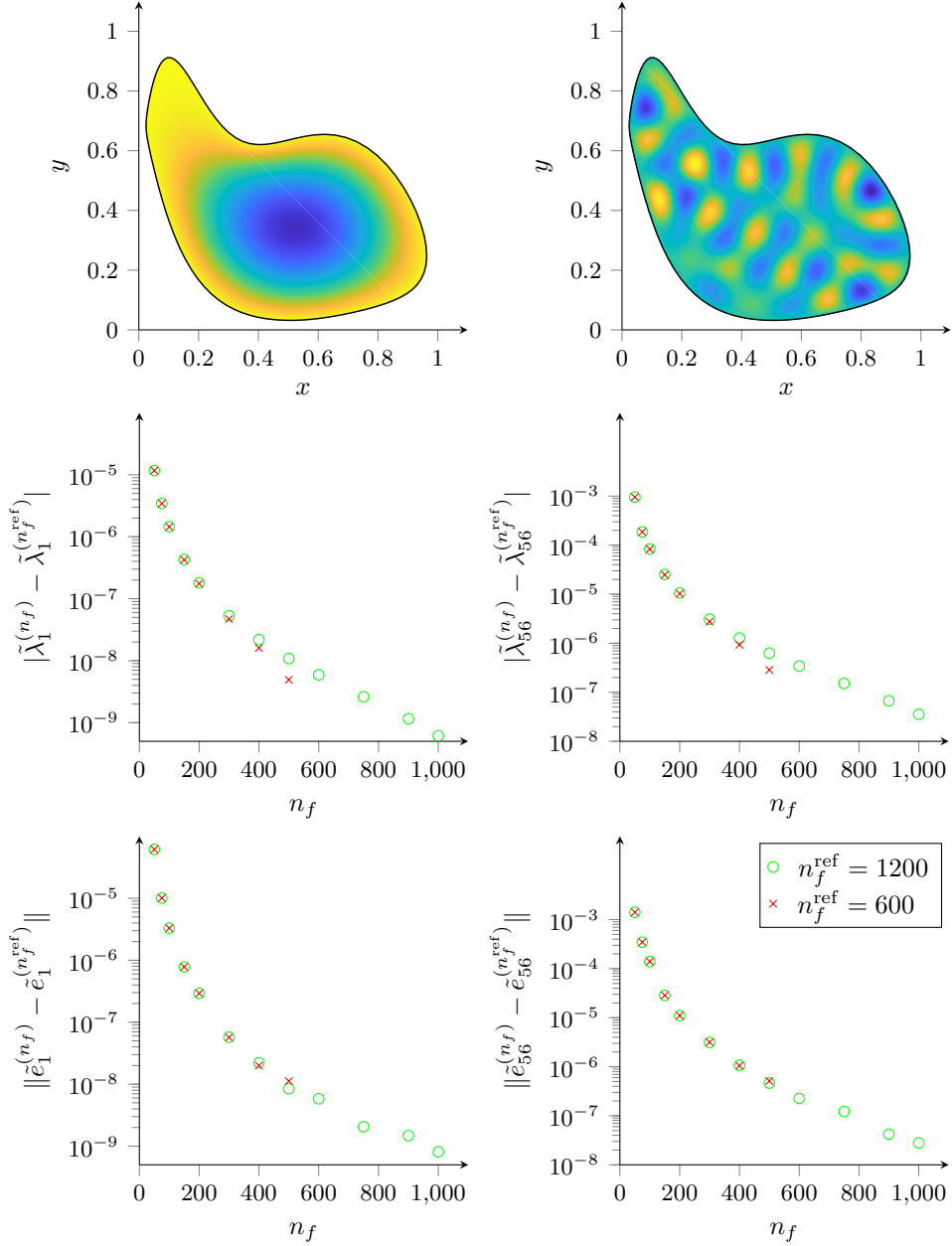
\begin{figure}[ht]
\newcommand{\Height}{5.98 cm}
\newcommand{\Width}{5.98 cm}
\newcommand{\PeanutWidth}{4.15 cm}
\newcommand{\X}{2.7}
\newcommand{\Yone}{1.0}
\newcommand{\Ytwo}{0.5}
\centering
 \begin{tabular}{@{}c@{}}
 % 1
 \hspace*{1.5em}{
 \begin{tikzpicture}
   \pgfdeclareimage[width=\PeanutWidth]{img}{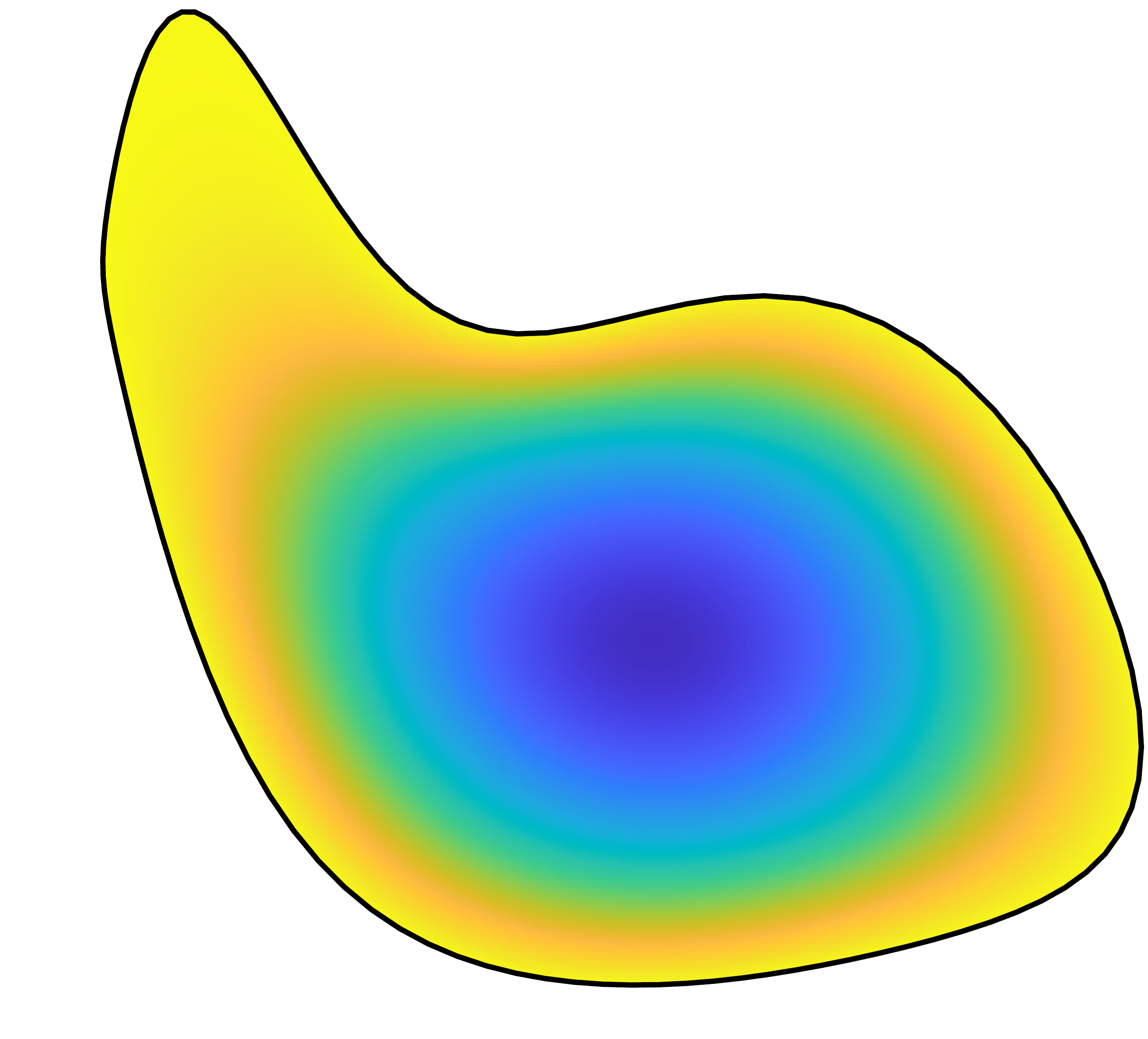}
        \node at (1.8,1.81) {\pgfuseimage{img}};
	\begin{axis}[
		height=\Height, 
		width=\Width,
		axis x line=bottom,
		axis y line=left,
		xlabel = {$x$},
		ylabel= {$y$},
		clip=false,
		ymin=0,ymax=1.1,
		xmin=0,xmax=1.1,
		xtick={0,0.2,0.4,0.6,0.8,1},
		ytick={0,0.2,0.4,0.6,0.8,1},
		tick label style={font=\footnotesize},
		tick align=outside
		]
	\end{axis}
\end{tikzpicture}
 }
 % 2
 \hspace*{1.5em}{
 \begin{tikzpicture}
   \pgfdeclareimage[width=\PeanutWidth]{img}{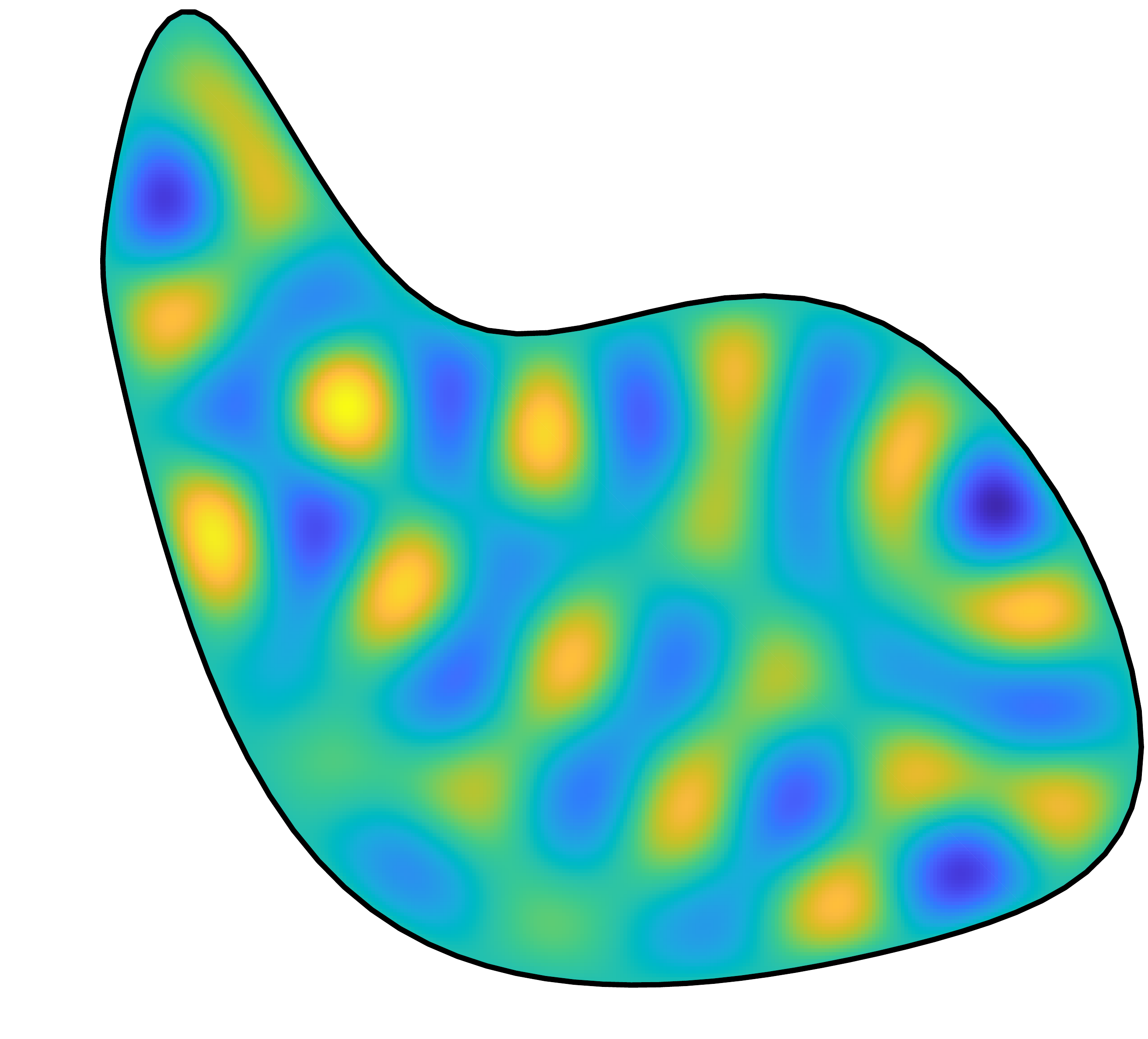}
        \node at (1.8,1.81) {\pgfuseimage{img}};
	\begin{axis}[
		height=\Height, 
		width=\Width,
		axis x line=bottom,
		axis y line=left,
		xlabel = {$x$},
		ylabel= {$y$},
		clip=false,
		ymin=0,ymax=1.1,
		xmin=0,xmax=1.1,
		xtick={0,0.2,0.4,0.6,0.8,1},
		ytick={0,0.2,0.4,0.6,0.8,1},
		tick label style={font=\footnotesize},
		tick align=outside
		]
	\end{axis}
\end{tikzpicture}
}
\\
 % 3
\begin{tikzpicture}
	\begin{axis}[
		height=\Height, 
		width=\Width,
		ymode = log,
		axis x line=bottom,
		axis y line=left,
		xlabel = {$n_f$},
		ylabel= {$|\tilde \lambda_1^{(n_f)} - \tilde \lambda_1^{(n_f^{\mathrm{ref}})} |$},
		clip=false,
		ymin=5e-10,ymax=1e-4,
		xmin=0,xmax=1100,
		xtick={0,200,400,600,800,1000},
		ytick={1e-9,1e-8,1e-7,1e-6,1e-5},
		tick label style={font=\footnotesize},
		tick align=outside
		]
		\addplot[color=green,only marks,mark=o] table {Peanut_ev_001_6.51554236_facemax=1000_reso=81.txt};
		\addplot[color=red,only marks,mark=x] table {Peanut_ev_001_6.51554236_facemax=500_reso=81.txt};
	\end{axis}
\end{tikzpicture}
% 4
\begin{tikzpicture}
	\begin{axis}[
		height=\Height, 
		width=\Width,
		ymode = log,
		axis x line=bottom,
		axis y line=left,
		xlabel = {$n_f$},
		ylabel= {$|\tilde \lambda_{56}^{(n_f)} - \tilde \lambda_{56}^{(n_f^{\mathrm{ref}})} |$},
		clip=false,
		ymin=1e-8,ymax=5e-2,
		xmin=0,xmax=1100,
		xtick={0,200,400,600,800,1000},
		ytick={1e-8,1e-7,1e-6,1e-5,1e-4,1e-3},
		tick label style={font=\footnotesize},
		tick align=outside
		]
		\addplot[color=green,only marks,mark=o] table {Peanut_ev_056_39.53663871_facemax=1000_reso=81.txt};
		\addplot[color=red,only marks,mark=x] table {Peanut_ev_056_39.53663872_facemax=500_reso=81.txt};
	\end{axis}
\end{tikzpicture}
\\
% 5
\begin{tikzpicture}
	\begin{axis}[
		height=\Height, 
		width=\Width,
		ymode = log,
		axis x line=bottom,
		axis y line=left,
		xlabel = {$n_f$},
		ylabel= {$\| \tilde e_1^{(n_f)} - \tilde e_1^{(n_f^{\mathrm{ref}})} \|$},
		clip=false,
		ymin=5e-10,ymax=1e-4,
		xmin=0,xmax=1100,
		xtick={0,200,400,600,800,1000},
		ytick={1e-9,1e-8,1e-7,1e-6,1e-5},
		tick label style={font=\footnotesize},
		tick align=outside
		]
		\addplot[color=green,only marks,mark=o] table {Peanut_L2_001_6.51554236_facemax=1000_reso=81.txt};
		\addplot[color=red,only marks,mark=x] table {Peanut_L2_001_6.51554236_facemax=500_reso=81.txt};
	\end{axis}
\end{tikzpicture}
% 6
\begin{tikzpicture}
	\begin{axis}[
		height=\Height, 
		width=\Width,
		ymode = log,
		axis x line=bottom,
		axis y line=left,
		xlabel = {$n_f$},
		ylabel= {$\| \tilde e_{56}^{(n_f)} - \tilde e_{56}^{(n_f^{\mathrm{ref}})} \|$},
		clip=false,
		ymin=1e-8,ymax=5e-2,
		xmin=0,xmax=1100,
		xtick={0,200,400,600,800,1000},
		ytick={1e-8,1e-7,1e-6,1e-5,1e-4,1e-3},
		tick label style={font=\footnotesize},
		tick align=outside
		]
		\addplot[color=green,only marks,mark=o] table {Peanut_L2_056_39.53663871_facemax=1000_reso=81.txt};
		\addplot[color=red,only marks,mark=x] table {Peanut_L2_056_39.53663872_facemax=500_reso=81.txt};
        \addlegendentry{$\;\; n_f^{\mathrm{ref}} = 1200$};
		\addlegendentry{$n_f^{\mathrm{ref}} = 600$};		
	\end{axis}
\end{tikzpicture}
\end{tabular}
\caption{Convergence results for two eigenfunctions of the Peanut shape (on the left, the first one, on the right, the 56th one), in red, the error is measured with respect to reference eigenvalues and eigenfunctions computed for $n_f^{\mathrm{ref}}=600$ boundary elements, in green, the reference data is computed for $n_f^{\mathrm{ref}}=1200$ boundary elements.}
\label{Reference1200600}
\end{figure}

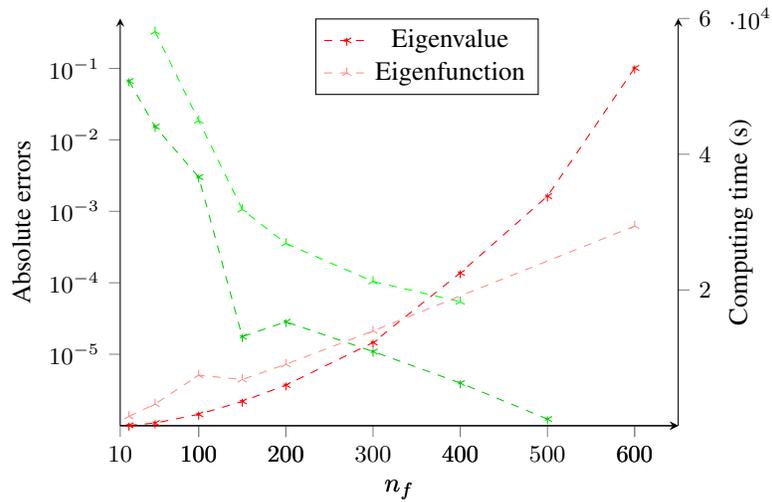
\begin{figure}
\newcommand{\Height}{7 cm}
\newcommand{\Width}{9 cm}
\newcommand{\X}{2.7}
\newcommand{\Yone}{1.0}
\newcommand{\Ytwo}{0.5}
\definecolor{dark_green}{RGB}{0,205,0}
\definecolor{light_green}{RGB}{0,255,0}
\definecolor{light_red}{RGB}{253,142,142}
\definecolor{dark_red}{RGB}{255,0,0}
\centering
 \begin{tabular}{@{}c@{}}
\begin{tikzpicture}
	\begin{axis}[
		height=\Height, 
		width=\Width,
		%xmode = log,
		ymode = log,
		axis x line=bottom,
		axis y line=left,
		xlabel = {$n_f$},
		ylabel= {Absolute errors},
		clip=false,
		ymin=1e-6,ymax=0.5,
		xmin=10,xmax=650,
		ytick={1e-1,1e-2,1e-3,1e-4,1e-5},
		tick label style={font=\footnotesize},
		tick align=outside
		]
		\addplot[color=dark_green,dashed,mark=asterisk] table {errorev_300.txt};
		\addplot[color=light_green,dashed,mark=Mercedes star] table {erroref_300.txt};	
	\end{axis}
	  	\begin{axis}[
		height=\Height, 
		width=\Width,
		%xmode = log,
		%ymode = log,
		axis x line=bottom,
		axis y line=right,
		xlabel = {$n_f$},
		ylabel= {Computing time (s)},
		clip=false,
		ymin=0,ymax=60000,
		xmin=10,xmax=650,
		xtick={10,100,200,400,600},
		ytick={20000,40000,60000},
		tick label style={font=\footnotesize},
		tick align=outside,
legend style={at={(0.75,1)}}] % <= to define position and font legend
% the following are the "images" and numbers in the legend
    %\addlegendimage{stealth-stealth,red,opacity=0.4} 
		]
		\addplot[color=dark_red,dashed,mark=asterisk] table 
		{compev_300.txt};
		\addplot[color=light_red,dashed,mark=Mercedes star] table {compef_300.txt};
		\addlegendentry{Eigenvalue};
		\addlegendentry{Eigenfunction};
	\end{axis}
\end{tikzpicture}
\end{tabular}
\caption{This plot shows, for different numbers of boundary elements $n_f$, absolute errors for the 300th eigenvalue and its corresponding eigenfunction in green (left axis) as well as CPU time for their computation in seconds (right axis). It can be seen that the CPU time for the eigenvalues behaves like $\mathcal{O}(n_f^2)$, while the time for the eigenfunctions is $\mathcal{O}(n_f)$. The computing times are very similar for other eigenvalues and eigenfunctions.}
\label{cputime_vs_errors}
\end{figure}

\begin{figure}
\newcommand{\Height}{7 cm}
\newcommand{\Width}{9 cm}
\newcommand{\X}{2.7}
\newcommand{\Yone}{1.0}
\newcommand{\Ytwo}{0.5}
\definecolor{dark_blue}{RGB}{0,0,200}
\definecolor{mid_blue}{RGB}{80,80,255}
\definecolor{light_blue}{RGB}{120,120,255}
\centering
 \begin{tabular}{@{}c@{}}
\hspace*{-3em}{\begin{tikzpicture}
	\begin{axis}[
		height=\Height, 
		width=\Width,
		%xmode = log,
		%ymode = log,
		axis x line=bottom,
		axis y line=left,
		xlabel = {$N$},
		ylabel= {Total CPU time (s)},
		clip=false,
		ymin=0,ymax=8e+7,
		xmin=0,xmax=550,
		ytick={2e+7,4e+7,6e+7},
		tick label style={font=\footnotesize},
		tick align=outside,
		legend pos = north west
		]
		\addplot[color=dark_blue,loosely dashed] table {totalcomptime_scenario_3.txt};
		\addplot[color=mid_blue,dashed] table {totalcomptime_scenario_2.txt};
		\addplot[color=light_blue,densely dashed] table {totalcomptime_scenario_1.txt};
		\addlegendentry{Scenario 3}
		\addlegendentry{Scenario 2}
		\addlegendentry{Scenario 1}
		\end{axis}
\end{tikzpicture}}
\end{tabular}
\caption{The total computation time on our machine for an orthonormal basis consisting of $N$ eigenvalues and eigenfunctions. Scenarios 1 to 3 correspond, from top to bottom, to the three pairs of error tolerances shown in Figure \ref{boundaryelements}.}
\label{totalbasecomptime}
\end{figure}
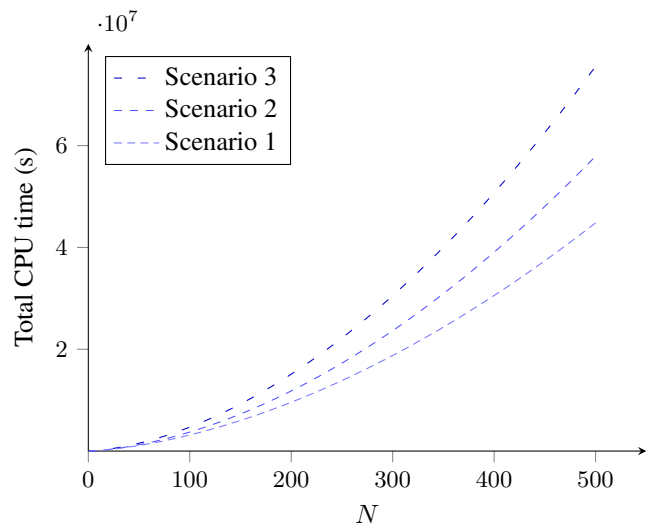

After precomputing the necessary numbers of boundary elements for each function (see Section \ref{eigenfunctionconvergence}) and the spectrum of $-\Delta$ on $D$ is known up to a few digits, the eigenvalues and eigenfunctions can be computed by applying the Beyn algorithm to each eigenvalue by enclosing it by a contour which does not enclose any other eigenvalue (see Figure \ref{contours}). The Peanut shape analyzed here is asymmetric and does not appear to have any multiple Dirichlet eigenvalues. If multiple eigenvalues occur, it can happen that the different eigenfunctions belonging to one eigenvalue are not returned by the algorithm as orthogonal functions, but as linear combinations of them.
This would necessitate an additional orthogonalization routine applied to those eigenfunctions, which would pose another error source.\\
Based on the convergence results for the Peanut shape, our base functions for the Peanut have been computed such that the error to a reference function with $n_f^{\mathrm{ref}}=600$ boundary elements stays below the tolerances $\epsilon^{(\lambda)}=2\cdot 10^{-4}$ and $\epsilon^{(\eta)} = 5\cdot 10^{-6}$.
In addition, we have computed some errors to even more accurate wavenumbers and eigenfunctions computed with $n_f^{\mathrm{ref}}=1,200$ boundary elements. A comparison of the reference errors is shown in Figure \ref{Reference1200600} for two wavenumbers. It can be seen that there is a slight discrepancy for the wavenumber error at $n_f=500$, but the reference error for the eigenfunctions is almost the same for both reference functions, so it can be assumed that the reference functions are sufficiently accurate.\\
On the computational cost, we note that for an orthonormal basis of a certain number $N$ of  eigenvalues and eigenfunctions the computational time with our implementation behaved asymptotically as $\mathcal{O}(N^2)$. In more detail, it is shown in Figure \ref{cputime_vs_errors} that, with growing number $n_f$ of used boundary elements, the computation times for an eigenvalue rises as $\mathcal{O}(n_f^2)$, and for an eigenfunction as $\mathcal{O}(n_f)$ (although it is more costly for low $n_f$). The computation time is very similar for different eigenvalues and the corresponding eigenfunctions, and neither increases nor decreases as the eigenvalue increases. Figure \ref{totalbasecomptime} shows an estimate of the total computational cost to compute an orthonormal basis of $N$ eigenvalues and eigenfunctions, i.e. the sum
\begin{align}\label{comptimesum}
    \sum_{i=1}^N \mathrm{ev}(n_f(\lambda_i)) + \mathrm{ef}(n_f(\lambda_i)),
\end{align}
where $n_f(\lambda)$ is the needed number of boundary elements for a given set of error tolerances, and $\mathrm{ev}(n_f)$ and $\mathrm{ef}(n_f)$ are the computation times for eigenvalue and eigenfunction using $n_f$ boundary elements. Based on empirical tests such as shown in Figure \ref{cputime_vs_errors}, we assumed that $\mathrm{ev}(n_f) = 0.15\cdot n_f^2$ and $\mathrm{ef}(n_f) = 50\cdot n_f$. Since $\lambda_i\sim i$, $n_f(i)\sim \sqrt{i}$ (see \ref{boundaryelements}), $\mathrm{ev}(i)\sim i^2$ and $\mathrm{ef}(i)\sim i$, the sum \eqref{comptimesum} contains a term behaving as $i$ and one behaving as $\sqrt{i}$, so that the whole sum is made up of two terms (one each for eigenvalues and eigenfunctions) behaving as $\mathcal{O}(N^2)+\mathcal{O}(N^{3/2})$.
\section{Error analysis.}\label{section_erroranalysis}

We present an error analysis for the strong error in the Galerkin exponential Euler scheme from \cite{KloedenRoyal} with approximated eigenfunctions and eigenvalues. We denote $\langle \cdot,\cdot \rangle := \langle \cdot,\cdot \rangle _{H} $ and $\| \cdot \| := \| \cdot \|_{H}$, and the exact solution at time $t$ by $U(t)$. Its approximation by the spectral Galerkin exponential Euler scheme is denoted by $V_k^{(1)}$. The approximation of $V_k^{(1)}$ by computing the inner products as numerical integrals over $ D$ is denoted by $V_k^{(2)}$. The approximation of $V_k^{(2)}$ using approximate eigenvalues is denoted by $V_k^{(3)}$. Finally, the approximation of $V_k^{(3)}$ using approximate eigenfunctions is denoted by $V_k^{(4)}$.
 We denote
\begin{align*}
 V_{k,j}^{(i)} := \langle V_k^{(i)}, e_j \rangle,\quad i=1,2,3,4,\quad 
 \epsilon_{k,j}^{(i)} := \mathbb{E}|V_{k,j}^{(i-1)}-V_{k,j}^{(i)}|,\quad i=2,3,4,
\end{align*}
with $k=0,\dots,M$ and $j=1,\dots,N$.
First, error estimates for $\|V_k^{(i)}-V_k^{(i+1)}\|$ are given for $i=1,2,3$ in Lemmas \ref{E2slemma}, \ref{E3slemma} and \ref{E4slemma} which are then combined to give a general result on the error $\|U(t_k)-V_k^{(4)}\|$ in Theorem \ref{strongerrortheorem}. Note that $V_k^{(4)}$ is the output of the computer programs, so the error estimate in Theorem \ref{strongerrortheorem} contains Fourier-Galerkin coefficients of $V_k^{(4)}$ which can then either be inserted for a run-specific bound, or be bound by an empirically found constant which is generally found to be an upper bound (which we have done in Example \ref{errorboundexample}).
 \begin{lemma}\label{E2slemma}
 We have
 \begin{align}\label{lemmabound2}
 \mathbb{E}\| V_k^{(1)}-V_k^{(2)} \| \leq \sum_{j=1}^N \epsilon_{k,j}^{(2)}
 \end{align}
 and, for a constant $C>0$, 
 \begin{align}\label{recursion1}
 \epsilon_{k+1,j}^{(2)} \leq \e^{-\lambda_j h} \epsilon_{k,j}^{(2)} + \frac{1-\e^{-\lambda_j h}}{\lambda_j}\cdot\left[L\cdot \sum_{i=1}^N \epsilon_{k,i}^{(2)} + C\cdot R^{-2}\right].
 \end{align}
 \end{lemma}
\begin{prf}
We compute, distinguishing $f_N^j(V_k^{(2)})$ computed as an integral and $\widetilde{f_N^j(V_k^{(2)})}$ computed as a numerical integral,
\begin{align}\label{wholefunctionerror}
\begin{split}
\mathbb{E}\| V_k^{(1)}-V_k^{(2)} \| &= \mathbb{E}\left \|\sum_{j=1}^N V_{k,j}^{(1)}\cdot e_j - \sum_{j=1}^N V_{k,j}^{(2)}\cdot e_j \right \| = \mathbb{E}\left\| \sum_{j=1}^N (V_{k,j}^{(1)}-V_{k,j}^{(2)})\cdot e_j \right\| \\
&\leq \sum_{j=1}^N \mathbb{E}|V_{k,j}^{(1)}-V_{k,j}^{(2)}| = \sum_{j=1}^N \epsilon_{k,j}^{(2)}
\end{split}
\end{align}
and
\begin{align*}
 \epsilon_{k+1,j}^{(2)}  &=\mathbb{E}\left| \e^{-\lambda_j h} V_{k,j}^{(1)} + \frac{1-\e^{-\lambda_j h}}{\lambda_j} f_N^j(V_k^{(1)}) + \left(\frac{q_j}{2\lambda_j}(1-\e^{-2\lambda_j h})\right)^{1/2}R_k^j\right. \\
 &\quad - \left.\left[ \e^{-\lambda_j h} V_{k,j}^{(2)} + \frac{1-\e^{-\lambda_j h}}{\lambda_j} \widetilde{f_N^j(V_k^{(2)})} + \left(\frac{q_j}{2\lambda_j}(1-\e^{-2\lambda_j h})\right)^{1/2}R_k^j\right] \right|\\
 &\leq\ \e^{-\lambda_j h} \epsilon_{k,j}^{(2)} + \frac{1-\e^{-\lambda_j h}}{\lambda_j}\mathbb{E}|f_N^j(V_k^{(1)})-\widetilde{f_N^j(V_k^{(2)})}|.
\end{align*}
%Let $\eta_k(x) := V_k^{(2)}(x)-V_k^{(1)}(x) = \sum_{j=1}^N \left[ V_{k,j}^{(2)}-V_{k,j}^{(1)} \right] e_j(x)$. Making use of the fact that in two dimensions, Simpson's numerical integration rule converges with order two,  we obtain
\begin{align*}
 &\mathbb{E}|f_N^j(V_k^{(1)})-\widetilde{f_N^j(V_k^{(2)})}| \leq \mathbb{E}|f_N^j(V_k^{(1)})-f_N^j(V_k^{(2)})| + \mathbb{E}|f_N^j(V_k^{(2)})-\widetilde{f_N^j(V_k^{(2)})}| \\ 
 \leq\ &\mathbb{E}\int_D \left[ f_N(V_k^{(1)}(x)) - f_N(V_k^{(1)}(x) + \eta_k(x))\right] e_j(x)\; \de x + C\cdot R^{-2}\\
 %\leq\ &\left[ \int_ D |f_N(V_k^{(1)}(x))-f_N(V_k^{(1)}(x)+\eta_k(x))|^2 e_j(x)\; \de x \right]^{1/2} + C\cdot R^{-2}\\
 \leq\ &\mathbb{E}\left[ \int_D C [L\cdot \eta_k(x)]^2\; \de x\right]^{1/2} + C\cdot R^{-2} = L\cdot \| \eta_k \| + C\cdot R^{-2} \\
 \leq\ &L\cdot \sum_{i=1}^N \epsilon_{k,i}^{(2)} + C\cdot R^{-2},
\end{align*}
where the H\"{o}lder inequality and the Lipschitz continuity of $f_N$ were used. In total, we have
\begin{equation*}
\pushQED{\qed} 
 \epsilon_{k+1,j}^{(2)} \leq \e^{-\lambda_j h} \epsilon_{k,j}^{(2)} + \frac{1-\e^{-\lambda_j h}}{\lambda_j} \cdot\left[L\cdot \sum_{i=1}^N \epsilon_{k,i}^{(2)} + C \cdot R^{-2}\right].
 \hfill\qed
\end{equation*}
\end{prf}
 \begin{lemma}\label{E3slemma}
 We have
 
 \begin{align}\label{lemmabound3}
  \mathbb{E}\| V_k^{(2)}-V_k^{(3)} \| \leq \sum_{j=1}^N \epsilon_{k,j}^{(3)}
 \end{align}
 and
 \begin{align}\label{recursion2}
 \epsilon_{k+1,j}^{(3)} &\leq \e^{-\lambda_j h}\epsilon_{k,j}^{(3)} + \e^{-\lambda_j h}|1-\e^{-\epsilon_j^{(\lambda)} h}|\cdot (\epsilon_{k,j}^{(4)} + \mathbb{E}|V_{k,j}^{(4)}|) \\
  &\quad + L\cdot \frac{1-\e^{-\lambda_j h}}{\lambda_j} \cdot \sum_{j=1}^N \epsilon_{k,j}^{(3)} \nonumber \\
  &\quad + \frac{1}{\lambda_j\tilde \lambda_j}\left[\epsilon^{(\lambda)} + |\e^{-\tilde \lambda_j h} \lambda_j - \e^{-\lambda_j h} \tilde \lambda_j |\right]\cdot \left[L\cdot \sum_{j=1}^N \epsilon_{k,j}^{(4)} + \mathbb{E}|f_N^j(V_k^{(4)})|  \right] \nonumber \\
 &\quad + \left| \sqrt{\frac{q_j}{2\lambda_j}(1-\e^{-2\lambda_j h})} - \sqrt{\frac{q_j}{2(\lambda_j+\epsilon_j^{(\lambda)})}(1-\e^{-2(\lambda_j+\epsilon_j^{(\lambda)}) h})}\right| \nonumber
 \end{align}
 \end{lemma}
 \begin{prf}
 Suppose that we have approximate eigenvalues $\tilde \lambda_j = \lambda_j+\epsilon_j^{(\lambda)}$, $j=1,\dots,N$.
Then by the same steps as in \eqref{wholefunctionerror}, the resulting error of the whole function is
\begin{align*}
 \mathbb{E}\| V_k^{(2)}-V_k^{(3)} \| \leq \sum_{j=1}^N \epsilon_{k,j}^{(3)}.
\end{align*}
The coefficient errors after the next time step are given by
\begin{align}
 \epsilon_{k+1,j}^{(3)} &= \mathbb{E}\left| \e^{-\lambda_j h}V_{k,j}^{(2)} + \frac{1-\e^{-\lambda_j h}}{\lambda_j}f_N^j(V_k^{(2)}) + \left(\frac{q_j}{2\lambda_j}(1-\e^{-\lambda_j h})\right)^{1/2} R_k^j \right. \nonumber\\
 &\quad - \left.\left[ \e^{-\tilde \lambda_j h}V_{k,j}^{(2)} + \frac{1-\e^{-\tilde \lambda_j h}}{\tilde \lambda_j}f_N^j(V_k^{(3)}) + \left(\frac{q_j}{2\tilde \lambda_j}(1-\e^{-2\tilde \lambda_j h})\right)^{1/2}R_k^j\right]\right| \nonumber\\
 &\leq \mathbb{E}|\e^{-\lambda_j h}V_{k,j}^{(2)} - \e^{-\tilde \lambda_j h} V_{k,j}^{(3)}| \nonumber\\
 &\quad + \mathbb{E}\left| \frac{1-\e^{-\lambda_j h}}{\lambda_j}f_N^j(V_k^{(2)}) - \frac{1-\e^{-\tilde \lambda_j h}}{\tilde \lambda_j}f_N^j(V_k^{(3)})\right|\nonumber \\
 &\quad + \mathbb{E}\left|\left[\left(\frac{q_j}{2\lambda_j}(1-\e^{-2\lambda_j h})\right)^{1/2} - \left(\frac{q_j}{2\tilde\lambda_j}(1-\e^{-2\tilde\lambda_j h})\right)^{1/2}\right]R_k^j\right|.\nonumber 
\end{align}
Going through this term by term, we have that
\begin{align}
&\quad \mathbb{E}|\e^{-\lambda_j h}V_{k,j}^{(2)} - \e^{-\tilde \lambda_j h} V_{k,j}^{(3)}| \nonumber\\
 &\leq \mathbb{E}|\e^{-\lambda_j h} V_{k,j}^{(2)}-\e^{-\lambda_j h}V_{k,j}^{(3)}| + \mathbb{E}|\e^{-\lambda_j h} V_{k,j}^{(3)}-\e^{-\tilde \lambda_j h}V_{k,j}^{(3)}| \nonumber\\
 &= \e^{-\lambda_j h} \mathbb{E}|V_{k,j}^{(2)}-V_{k,j}^{(3)}| + |\e^{-\lambda_j h}-\e^{-(\lambda_j+\epsilon_j^{(\lambda)}) h}|\cdot \mathbb{E}|V_{k,j}^{(3)}| \nonumber\\
 &\leq \e^{-\lambda_j h}\epsilon_{k,j}^{(3)} + \e^{-\lambda_j h}|1-\e^{-\epsilon_j^{(\lambda)} h}|\cdot \mathbb{E}|V_{k,j}^{(3)}|\nonumber\\
 &\leq \e^{-\lambda_j h}\epsilon_{k,j}^{(3)} + \e^{-\lambda_j h}|1-\e^{-\epsilon_j^{(\lambda)} h}|\cdot (\mathbb{E}|V_{k,j}^{(3)}-V_{k,j}^{(4)}| + \mathbb{E}|V_{k,j}^{(4)}|)\nonumber\\
 &\leq \e^{-\lambda_j h}\epsilon_{k,j}^{(3)} + \e^{-\lambda_j h}|1-\e^{-\epsilon_j^{(\lambda)} h}|\cdot (\epsilon_{k,j}^{(4)} + \mathbb{E}|V_{k,j}^{(4)}|),\label{E31s}
\end{align}
where we have made the last estimate since unlike $\mathbb{E}[V_k^{(3)}]$, $\mathbb{E}[V_k^{(4)}]$ can be approximated from the numerical implementation. Next, we have
\begin{align}
&\quad \mathbb{E}\left| \frac{1-\e^{-\lambda_j h}}{\lambda_j}f_N^j(V_k^{(2)}) - \frac{1-\e^{-\tilde \lambda_j h}}{\tilde \lambda_j}f_N^j(V_k^{(3)})\right| \nonumber\\
 &\leq \mathbb{E}\left|\frac{1-\e^{-\lambda_j h}}{\lambda_j}f_N^j(V_k^{(2)})-\frac{1-\e^{-\lambda_j h}}{\lambda_j}f_N^j(V_k^{(3)})\right| \nonumber\\
 &\quad + \mathbb{E}\left| \frac{1-\e^{-\lambda_j h}}{\lambda_j}f_N^j(V_k^{(3)})-\frac{1-\e^{-\tilde\lambda_j h}}{\tilde \lambda_j}f_N^j(V_k^{(3)})\right| \nonumber\\
 &= \frac{1-\e^{-\lambda_j h}}{\lambda_j}\mathbb{E}|f_N^j(V_k^{(2)})-f_N^j(V_k^{(3)})| + \left|\frac{1-\e^{-\lambda_j h}}{\lambda_j}-\frac{1-\e^{-\tilde \lambda_j h}}{\tilde \lambda_j}\right|\cdot \mathbb{E}|f_N^j(V_k^{(3)})| \nonumber\\
 &\leq  L\cdot \frac{1-\e^{-\lambda_j h}}{\lambda_j} \cdot \sum_{j=1}^N \epsilon_{k,j}^{(3)} \nonumber\\
 &\quad + \frac{1}{\lambda_j\tilde \lambda_j}\left[\epsilon^{(\lambda)} + |\e^{-\tilde \lambda_j h} \lambda_j - \e^{-\lambda_j h} \tilde \lambda_j |\right]\cdot \left[\mathbb{E}|f_N^j(V_k^{(3)})-f_N^j(V_k^{(4)})| + \mathbb{E}|f_N^j(V_k^{(4)})|  \right] \\
 \label{E32s}
 &\leq L\cdot \frac{1-\e^{-\lambda_j h}}{\lambda_j} \cdot \sum_{j=1}^N \epsilon_{k,j}^{(3)} \\
 &\quad + \frac{1}{\lambda_j\tilde \lambda_j}\left[\epsilon^{(\lambda)} + |\e^{-\tilde \lambda_j h} \lambda_j - \e^{-\lambda_j h} \tilde \lambda_j |\right]\cdot \left[L\cdot \sum_{j=1}^N \epsilon_{k,j}^{(4)} + \mathbb{E}|f_N^j(V_k^{(4)})|  \right].\nonumber
\end{align}
Finally, we have (using that $\mathbb{E}|X| \leq [\mathbb{E}|X|^2]^{1/2}$)
\begin{align}
&\mathbb{E}\left|\left[\left(\frac{q_j}{2\lambda_j}(1-\e^{-2\lambda_j h})\right)^{1/2} - \left(\frac{q_j}{2\tilde\lambda_j}(1-\e^{-2\tilde\lambda_j h})\right)^{1/2}\right]R_k^j\right| \nonumber\\
\label{E33s}
   &=\left| \sqrt{\frac{q_j}{2\lambda_j}(1-\e^{-2\lambda_j h})} - \sqrt{\frac{q_j}{2(\lambda_j+\epsilon_j^{(\lambda)})}(1-\e^{-2(\lambda_j+\epsilon_j^{(\lambda)}) h})}\right| \mathbb{E}|(R_k^j)^2|^{1/2},
\end{align}
 and $\mathbb{E}|(R_k^j)^2| = \text{Var}(R_k^j) = 1$. Summing \eqref{E31s}, \eqref{E32s} and \eqref{E33s}, we obtain \eqref{recursion2}.\qed
 \end{prf}
 \begin{lemma}\label{E4slemma}
 We have 
 
 \begin{align}\label{lemmabound4}
 \mathbb{E}\| V_k^{(3)}-V_k^{(4)}\|\leq \sum_{j=1}^N\left[ \epsilon_{k,j}^{(4)} + \epsilon^{(\eta)}\mathbb{E}|V_{k,j}^{(4)}|\right]
 \end{align}
 and
\begin{align}\label{recursion3}
 \begin{split}
 \epsilon_{k+1,j}^{(4)} &\leq \e^{-\tilde\lambda_j h}\epsilon_{k,j}^{(4)} + \e^{-\tilde\lambda_j h}\epsilon_j^{(\eta)}\mathbb{E}\|V_k^{(4)}\| \\
 &\quad + \frac{1-\e^{-\tilde\lambda_j h}}{\tilde \lambda_j}\left[L\cdot \sum_{j=1}^N \epsilon_{k,j}^{(4)} + \epsilon_j^{(\eta)}\cdot \mathbb{E}\|f_N(V_k^{(4)})\| \right].
 \end{split}
 \end{align}
 \end{lemma}
 \begin{prf}
 Suppose that we have approximate eigenfunctions
\begin{align*}
 \tilde e_j = e_j+\eta_j,\quad \|\eta_j\|\leq \epsilon_j^{(\eta)}.
\end{align*}
Then, the Fourier-Galerkin expansion of $V_k^{(4)}$ also involves the approximate eigenfunctions, i.e. $V_k^{(4)} = \sum_{j=1}^N (e_j+\eta_j)V_{k,j}^{(4)}$, and the error for the function is 
\begin{alignat*}{2}
 & \mathbb{E}\| V_k^{(3)}-V_k^{(4)}\| \\
 &= \mathbb{E}\left\| \sum_{j=1}^N e_j V_{k,j}^{(3)} - \sum_{j=1}^N (e_j+\eta_j)V_{k,j}^{(4)} \right\| &&= \mathbb{E}\left \| \sum_{j=1}^N e_j (V_{k,j}^{(3)}-V_{k,j}^{(4)}) - \sum_{j=1}^N \eta_j V_{k,j}^{(4)} \right\| \\
 &\leq \;\;\sum_{j=1}^N \epsilon_{k,j}^{(4)} + \sum_{j=1}^N \epsilon^{(\eta)}\cdot \mathbb{E}|V_{k,j}^{(4)}| &&=\;\; \sum_{j=1}^N\left[ \epsilon_{k,j}^{(4)} + \epsilon^{(\eta)}\mathbb{E}|V_{k,j}^{(4)}|\right].
\end{alignat*}
The error for the next coefficient errors is (with the wide tilde denoting inner products with $e_j+\eta_j$ instead of with $e_j$)
\begin{align*}
 \epsilon_{k+1,j}^{(4)} &= \mathbb{E}\left| \e^{-\tilde\lambda_j h} V_{k,j}^{(3)} + \frac{1-\e^{-\tilde\lambda_j h}}{\tilde \lambda_j}f_N^j(V_k^{(3)}) + \left(\frac{q_j}{2\tilde\lambda_j}(1-\e^{-2\tilde\lambda_j h})\right)^{1/2}R_k^j \right.\\
&\quad - \left.\left[ \e^{-\tilde\lambda_j h} \widetilde{V_{k,j}^{(4)}} + \frac{1-\e^{-\tilde\lambda_j h}}{\tilde \lambda_j}\widetilde{f_N^j(V_k^{(4)})} + \left(\frac{q_j}{2\tilde\lambda_j}(1-\e^{-2\tilde\lambda_j h})\right)^{1/2}R_k^j\right]\right| \\
&= \mathbb{E}\left|\e^{-\tilde \lambda_j h}(V_{k,j}^{(3)}-\widetilde{V_{k,j}^{(4)}}) + \frac{1-\e^{-\tilde \lambda_j h}}{\tilde \lambda_j}[f_N^j(V_k^{(3)})-\widetilde{f_N^j(V_k^{(4)})}]\right|.
\end{align*}
 Dealing with these terms separately, we have
 \begin{align}
   \ &\e^{-\tilde\lambda_j h}\mathbb{E}|V_{k,j}^{(3)}-V_{k,j}^{(4)}| + \e^{-\tilde\lambda_j h}\mathbb{E}|V_{k,j}^{(4)}-\widetilde{V_{k,j}^{(4)}}| \nonumber\\
  \leq\ &\e^{-\tilde\lambda_j h}\epsilon_{k,j}^{(4)} + \e^{-\tilde\lambda_j h}\mathbb{E}|\langle V_k^{(4)},e_j \rangle - \langle V_k^{(4)},e_j+\eta_j\rangle | \nonumber\\
  =\ &\e^{-\tilde\lambda_j h}\epsilon_{k,j}^{(4)} + \e^{-\tilde\lambda_j h}\mathbb{E}|\langle V_k^{(4)},-\eta_j \rangle | \nonumber\\
  \label{E41s}
  \leq\ &\e^{-\tilde\lambda_j h}\epsilon_{k,j}^{(4)} + \e^{-\tilde\lambda_j h}\epsilon_j^{(\eta)}\mathbb{E}\|V_k^{(4)}\|
 \end{align}
and
\begin{align}
 &\frac{1-\e^{-\tilde\lambda_j h}}{\tilde \lambda_j}\left[\mathbb{E}| f_N^j(V_k^{(3)})-f_N^j(V_k^{(4)})| + \mathbb{E}|f_N^j(V_k^{(4)})-\widetilde{f_N^j(V_k^{(4)})}| \right] \nonumber\\
 \leq\ &\frac{1-\e^{-\tilde\lambda_j h}}{\tilde \lambda_j}\left[\mathbb{E}| f_N^j(V_k^{(3)})-f_N^j(V_k^{(4)})| + \mathbb{E}|\langle f_N(V_k^{(4)}),e_j\rangle - \langle f_N(V_k^{(4)}),e_j+\eta_j\rangle| \right]\nonumber\\
 \leq\ &\frac{1-\e^{-\tilde\lambda_j h}}{\tilde \lambda_j}\left[\mathbb{E}| f_N^j(V_k^{(3)})-f_N^j(V_k^{(4)})| + \mathbb{E}|\langle f_N(V_k^{(4)}),-\eta_j\rangle| \right]\nonumber\\
 \label{E42s}
 \leq\ &\frac{1-\e^{-\tilde\lambda_j h}}{\tilde \lambda_j}\left[L\cdot \sum_{j=1}^N \epsilon_{k,j}^{(4)} + \epsilon_j^{(\eta)}\cdot \mathbb{E}\|f_N(V_k^{(4)})\| \right].
\end{align}
Summing \eqref{E41s} and \eqref{E42s}, we obtain \eqref{recursion3}.\qed
 \end{prf}
 
 \begin{theorem}\label{strongerrortheorem}
  Suppose that we have fixed parameters $T$, $N$ eigenfunctions, $M$ time steps, a grid resolution $R$, a Lipschitz constant $L=\| f' \|_\infty$, Dirichlet eigenvalues $0<\lambda_1<\lambda_2\leq \lambda_3\leq \dots \leq \lambda_n$, 
eigenvalue error tolerance $\epsilon^{(\lambda)}$, and eigenfunction error tolerance $\epsilon^{(\eta)}$. We denote the Fourier coefficients $V_{k,j}^{(4)}$ of the approximated solution $V^{(4)}_k$ and the Fourier coefficients $f_N^j(V_k^{(4)})$ of $f_N(U_{t_k})$ in the matrices $\bar{V} = (V_{k,j}^{(4)})$, $f_N(\bar{V}) = (f_N^j(V_{k}^{(4)}))$, $k=1,\dots,M$, $j=1,\dots,N$.
Then, for any $\epsilon\in (0,1)$, we have a bound for the strong error of the complete scheme at time $T$
\begin{align}
\begin{split}
\mathbb{E}\|U(T)-V_M^{(4)}\| \leq \mathcal{E}^s(T,N,M,R,L,\lambda_1,\lambda_N,\epsilon^{(\lambda)},\epsilon^{(\eta)},\bar{V},f_N(\bar{V})) = \\ 
\label{errorbound}
\mathcal{E}_M^{(1)}(T,N,M,\lambda_N) + \mathcal{E}_M^{(234)}(T,N,M,R,L,\lambda_1,\lambda_N,\epsilon^{(\lambda)},\epsilon^{(\eta)},\bar{V},f_N(\bar{V})),
\end{split}
\end{align}
where, for a constant $C_T>0$, 
\begin{align}\label{errorjentzenpart}
\mathcal{E}_M^{(1)}(T,N,M,\lambda_N)\leq C_T\left(\lambda_N^{\epsilon-1}+\frac{\log(M)}{M}\right) 
\end{align}
and
\begin{align}
&\mathcal{E}_M^{(234)}(T,N,M,R,L,\lambda_1,\epsilon^{(\lambda)},\epsilon^{(\eta)},\bar{V},f_N(\bar{V}))\nonumber \\
\label{bound1}
&\quad \leq a_2^M\cdot \epsilon_0^{(2)} + b_2\cdot \frac{1-a_2^M}{1-a_2} + a_2^M\cdot \epsilon_0^{(3)} + b_3 \cdot \frac{1-a_2^M}{1-a_2} + \tilde a_2^M\cdot \epsilon_0^{(4)} + b_4\cdot \frac{1-\tilde a_2^M}{1-\tilde a_2} + b_5 \\ 
&\quad \rightarrow \frac{b_2}{1-a_2} + \frac{b_3}{1-a_2} + \frac{b_4}{1-\tilde a_2} + b_5\ \text{as}\ M\rightarrow\infty,\ \text{if}\ a_2<1\ (\text{since }\tilde a_2\leq a_2).\nonumber
\end{align}
where, for some constant $C>0$,
\begin{align*}
a_2 &= \e^{-\lambda_1 h} + L\cdot N\cdot \frac{1-\e^{-\lambda_1 h}}{\lambda_1},\quad \tilde a_2 = \e^{-\tilde \lambda_1 h} + L\cdot N\cdot \frac{1-\e^{-\tilde \lambda_1 h}}{\tilde \lambda_1}, \nonumber\\
b_2 &= \frac{1-\e^{-\lambda_1 h}}{\lambda_1}\cdot N\cdot C\cdot R^{-2},\\
b_3 &= \e^{-\lambda_1 h} |1-\e^{-\epsilon^{(\lambda)} h} |\cdot \sup_{k=0,\dots,M}\left[ \sum_{j=1}^N (\epsilon_{k,j}^{(4)} + \mathbb{E}|V_{k,j}^{(4)}|) \right] \\
&\quad + \frac{1}{\lambda_1 \tilde \lambda_1}\left[ \epsilon^{(\lambda)} + |\e^{-\tilde \lambda_1 h} \lambda_1 - \e^{-\lambda_1 h} \tilde \lambda_1 | \right] \\
&\quad \cdot \sup_{k=0,\dots,M} \left[N\cdot L\cdot \sum_{j=1}^N \epsilon_{k,j}^{(4)} + \sum_{j=1}^N \mathbb{E}|f_N^j(V_k^{(4)})| \right] \\ 
&\quad + N\cdot \left| \sqrt{\frac{1}{2\lambda_1}(1-\e^{-2\lambda_1 h})}-\sqrt{\frac{1}{2(\lambda_1+\epsilon^{(\lambda)})}(1-\e^{-2(\lambda_1+\epsilon^{(\lambda)} )h})}\right| \\
b_4 &= N \cdot \sup_{j=1,\dots,N}\Biggl\{ (1+\mathbb{E}|V^{(4)}_{M,j}|) \\ 
 &\quad \cdot \left[ \e^{-\tilde \lambda_1 h}\cdot \epsilon^{(\eta)}\cdot \sup_{k=0,\dots,M} \mathbb{E}\| V_k^{(4)} \| + \epsilon^{(\eta)}\cdot \frac{1-\e^{-\tilde \lambda_1 h}}{\tilde \lambda_1} \cdot \sup_{k=0,\dots,M} \mathbb{E}\| f_N(V_k^{(4)}) \|  \right]\Biggl \},\\
 b_5 &= \epsilon^{(\eta)}\cdot \sum_{j=1}^N \mathbb{E}|V_{M,j}^{(4)}|.
\end{align*}
\end{theorem}
\begin{remark}
 For the case of the simpler error bound with $M\rightarrow \infty$, the condition $a_2<1$ is equivalent to $L\cdot N<\lambda_1$. This turns out to be a rather strict condition (see Section \ref{numericalexperiments}). The error bound \eqref{errorbound} requires very strict conditions on the parameters in order to be usable, we address this in Example \ref{errorboundexample}. This strictness stems from some crude estimates in the proofs to simplify expressions, for example the step of replacing all $\lambda_n$ by $\lambda_1$ in the proof of Theorem \ref{strongerrortheorem}.  
 \end{remark}
\begin{prf}
The error of the scheme can be split up as
\begin{align*}
 & \mathbb{E}(\|U(t_k)-V_k^{(4)}\|) \\
 &\leq \mathbb{E}(\|U(t_k)-V_k^{(1)}\|) + \mathbb{E}(\|V_k^{(1)}-V_k^{(2)}\|) + \mathbb{E}(\|V_k^{(2)}-V_k^{(3)}\|) + \mathbb{E}(\|V_k^{(3)}-V_k^{(4)}\|)\\
 &=: \mathbb{E}(\|U(t_k)-V_k^{(1)}\|) + \epsilon_k^{(2)} + \epsilon_k^{(3)} + \epsilon_k^{(4)} \leq \mathcal{E}_k^{(1)} + \mathcal{E}_k^{(2)} + \mathcal{E}_k^{(3)} + \mathcal{E}_k^{(4)},\quad k=1,\dots, M,
 \end{align*}
where the error bound $\mathcal{E}_M^{(1)}$ in \eqref{errorjentzenpart} is from \cite{KloedenRoyal}. The other terms relate to the second error term above as
\begin{align*}
\mathcal{E}_k^{(234)} = \mathcal{E}_k^{(2)} + \mathcal{E}_k^{(3)} + \mathcal{E}_k^{(4)},\quad k=1,\dots,M.
\end{align*}
The error bounds $\mathcal{E}_k^{(i)}$, $i=2,3,4$, are established in Lemmas \ref{E2slemma}, \ref{E3slemma} and \ref{E4slemma} as
\begin{align*}
 \epsilon_k^{(2)} \leq \sum_{j=1}^N \epsilon_{k,j}^{(2)} =: \mathcal{E}_k^{(2)},\quad
 \epsilon_k^{(3)}  \leq \sum_{j=1}^N \epsilon_{k,j}^{(3)} =: \mathcal{E}_k^{(3)},\quad 
 \epsilon_k^{(4)}  \leq \sum_{j=1}^N (1+|V_{k,j}^{(4)}|) \epsilon_{k,j}^{(4)} =: \mathcal{E}_k^{(4)}.
\end{align*}
Furthermore, the inequalities \eqref{recursion1}, \eqref{recursion2} and \eqref{recursion3} establish relations between $\epsilon_{k,j}^{(i)}$ and $\epsilon_{k+1,j}^{(i)}$ for $k=0,\dots,M-1$, $j=1,\dots,N$ and $i=2,3,4$. First, we consider in each time step the sum of all errors (summing over index $j$) and their bounds given by
\begin{alignat*}{2}
\sum_{j=1}^N \epsilon_{k+1,j}^{(2)} &\leq\ \sum_{j=1}^N\e^{-\lambda_j h} \epsilon_{k,j}^{(2)} + \sum_{j=1}^N \frac{1-\e^{-\lambda_j h}}{\lambda_j}\left[ L\cdot \sum_{i=1}^N \epsilon_{k,i}^{(2)} + D\cdot R^{-2}\right],\\
\sum_{j=1}^N \epsilon_{k+1,j}^{(3)} &\leq\ \sum_{j=1}^N \e^{-\lambda_j h}\epsilon_{k,j}^{(3)} + \sum_{j=1}^N \e^{-\lambda_j h}|1-\e^{-\epsilon^{(\lambda)} h}|\cdot (\epsilon_{k,j}^{(4)} + \mathbb{E}|V_{k,j}^{(4)}|) \\ 
&\quad + \sum_{j=1}^N L\cdot \frac{1-\e^{-\lambda_j h}}{\lambda_j} \cdot \sum_{i=1}^N \epsilon_{k,i}^{(3)} \\
 &\quad + \sum_{j=1}^N\frac{1}{\lambda_j\tilde \lambda_j}\left[\epsilon^{(\lambda)} + |\e^{-\tilde \lambda_j h} \lambda_j - \e^{-\lambda_j h} \tilde \lambda_j |\right]\cdot \left[L\cdot \sum_{i=1}^N \epsilon_{k,i}^{(4)} + \mathbb{E}|f_N^j(V_k^{(4)})|  \right] \\
&\quad + \sum_{j=1}^N\left| \sqrt{\frac{q_j}{2\lambda_j}(1-\e^{-2\lambda_j h})} - \sqrt{\frac{q_j}{2(\lambda_j+\epsilon^{(\lambda)})}(1-\e^{-2(\lambda_j+\epsilon^{(\lambda)}) h})}\right|,\\
\sum_{j=1}^N\epsilon_{k+1,j}^{(4)} \leq\ &\sum_{j=1}^N\e^{-\tilde\lambda_j h}\epsilon_{k,j}^{(4)} +\sum_{j=1}^N \e^{-\tilde\lambda_j h}\epsilon^{(\eta)}\mathbb{E}\|V_k^{(4)}\| \\
&\quad + \sum_{j=1}^N\frac{1-\e^{-\tilde\lambda_j h}}{\tilde \lambda_j}\left[L\cdot \sum_{i=1}^N \epsilon_{k,i}^{(4)} + \epsilon^{(\eta)}\cdot \mathbb{E}\|f_N(V_k^{(4)})\| \right].
\end{alignat*}
Since the terms in the above expressions involving $\lambda_j$ like $\e^{-\lambda_j h}$, $\frac{1-\e^{-\lambda_j h}}{\lambda_j}$,\linebreak $\frac{1}{\lambda_j\tilde \lambda_j}\left[\epsilon^{(\lambda)} + |\e^{-\tilde \lambda_j h} \lambda_j - \e^{-\lambda_j h} \tilde \lambda_j |\right]$ and 
\begin{equation*}
\textstyle
\left| \sqrt{\frac{q_j}{2\lambda_j}(1-\e^{-2\lambda_j h})} - \sqrt{\frac{q_j}{2(\lambda_j+\epsilon^{(\lambda)})}(1-\e^{-2(\lambda_j+\epsilon^{(\lambda)}) h})}\right|
\end{equation*}
are monotone decreasing in $\lambda_j$, we can bound all the above terms from above by replacing all $\lambda_j$ by $\lambda_1$. We obtain
\begin{align}
 \sum_{j=1}^N \epsilon_{k+1,j}^{(2)} &\leq \e^{-\lambda_1 h} \sum_{j=1}^N \epsilon_{k,j}^{(2)} + N\cdot L \cdot \frac{1-\e^{-\lambda_1 h}}{\lambda_1} \sum_{j=1}^N \epsilon_{k,j}^{(2)} + \frac{1-\e^{-\lambda_1 h}}{\lambda_1}\cdot N\cdot D \cdot R^{-2} \nonumber \\
 \label{iteration1}
 &= \left[ \e^{-\lambda_1 h} + N\cdot L\cdot \frac{1-\e^{-\lambda_1 h}}{\lambda_1}\right]\cdot \sum_{j=1}^N \epsilon_{k,j}^{(2)} + \frac{1-\e^{-\lambda_1 h}}{\lambda_1}\cdot N\cdot D\cdot R^{-2},\\
 \sum_{j=1}^N \epsilon_{k+1,j}^{(3)} &\leq \e^{-\lambda_1 h} \sum_{j=1}^N \epsilon_{k,j}^{(3)} + \e^{-\lambda_1 h} |1-\e^{-\epsilon^{(\lambda)} h}| \sum_{j=1}^N(\epsilon_{k,j}^{(4)} + \mathbb{E}|V_{k,j}^{(4)}|)\nonumber \\
 &\quad + N\cdot L\cdot \frac{1-\e^{-\lambda_1 h}}{\lambda_1}\cdot \sum_{j=1}^N\epsilon_{k,j}^{(3)} \nonumber\\
 &\quad + \frac{1}{\lambda_1 \tilde \lambda_1} \\
 &\quad \cdot \left[ \epsilon^{(\lambda)} + |\e^{-\tilde \lambda_1 h}\lambda_1 - \e^{-\lambda_1 h} \tilde \lambda_1 |\right]\cdot \left[N\cdot L\cdot \sum_{j=1}^N\epsilon_{k,j}^{(4)} + \sum_{j=1}^N\mathbb{E}|f_N^j(V_k^{(4)})| \right]\nonumber \\
  &\quad + N\cdot \left| \sqrt{\frac{1}{2\lambda_1}(1-\e^{-2\lambda_1 h})}-\sqrt{\frac{1}{2(\lambda_1+\epsilon^{(\lambda)})}(1-\e^{-2(\lambda_1+\epsilon^{(\lambda)} h)})}\right| \nonumber\\
   &\leq \left[ \e^{-\lambda_1 h} + N\cdot L\cdot \frac{1-\e^{-\lambda_1 h}}{\lambda_1}\right]\cdot \sum_{j=1}^N \epsilon_{k,j}^{(3)}   \label{iteration2}\\
 &\quad + \e^{-\lambda_1 h} |1-\e^{-\epsilon^{(\lambda)} h} |\cdot \sup_{k=0,\dots,M}\left[ \sum_{j=1}^N (\epsilon_{k,j}^{(4)} + \mathbb{E}|V_{k,j}^{(4)}|) \right] \nonumber\\
 &\quad + \frac{1}{\lambda_1 \tilde \lambda_1}\left[ \epsilon^{(\lambda)} + |\e^{-\tilde \lambda_1 h} \lambda_1 - \e^{-\lambda_1 h} \tilde \lambda_1 | \right]\nonumber\\
 &\quad \cdot \sup_{k=0,\dots,M} \left[N\cdot L\cdot \sum_{j=1}^N \epsilon_{k,j}^{(4)} + \sum_{j=1}^N \mathbb{E}|f_N^j(V_k^{(4)})| \right]\nonumber \\
 &\quad + N\cdot \left| \sqrt{\frac{1}{2\lambda_1}(1-\e^{-2\lambda_1 h})}-\sqrt{\frac{1}{2(\lambda_1+\epsilon^{(\lambda)})}(1-\e^{-2(\lambda_1+\epsilon^{(\lambda)} h)})}\right|,\nonumber
 \end{align}
 \begin{align}
 \sum_{j=1}^N \epsilon_{k+1,j}^{(4)} &\leq \e^{-\tilde \lambda_1 h} \sum_{j=1}^N \epsilon_{k,j}^{(4)} + \e^{-\tilde \lambda_1 h} \cdot N\cdot \epsilon^{(\eta)} \cdot \mathbb{E}\| V_k^{(4)} \| \nonumber\\
 &\quad + N\cdot \frac{1-\e^{-\tilde \lambda_1 h}}{\tilde \lambda_1} \left[ L\cdot \sum_{j=1}^N \epsilon_{k,j}^{(4)} + \epsilon^{(\eta)} \cdot \mathbb{E}\| f_N(V_k^{(4)}) \|\right] \nonumber\\
 &\leq \left[ \e^{-\tilde \lambda_1 h} + N\cdot L\cdot \frac{1-\e^{-\tilde \lambda_1 h}}{\tilde \lambda_1}\right]\cdot \sum_{j=1}^N \epsilon_{k,j}^{(4)} \label{iteration3}\\
 &\quad + \e^{-\tilde \lambda_1 h}\cdot N \cdot \epsilon^{(\eta)}\cdot \sup_{k=0,\dots,M} \mathbb{E}\| V_k^{(4)} \| \\
 &\quad + N\cdot \epsilon^{(\eta)} \cdot \frac{1-\e^{-\tilde \lambda_1 h}}{\tilde \lambda_1} \cdot \sup_{k=0,\dots,M} \mathbb{E}\| f_N(V_k^{(4)}) \|.\nonumber
\end{align}
We can now see that the error sum at $t_{k+1}$ can be bounded by an expression which is the image of the affine map
\begin{align*}
\psi_{a,b}: \R\rightarrow \R,\quad \psi_{a,b}(x) = ax+b
\end{align*}
for $a>0,\ b\in\R$. For a given start value, the behavior of the iteration $\psi_{a,b}^M(x_0)$, $M>0$, depends on whether $a<1$, $a=1$ or $a>1$. We have the following cases:
\begin{numcases}{\psi_{a,b}^M(x_0) = \ }\label{limitingbehaviorstart}
 a^Mx_0 + b\sum_{i=1}^{M-1} a^i = a^Mx_0 + b\cdot \frac{1-a^M}{1-a} \stackrel{M\rightarrow \infty}{\longrightarrow}  \frac{b}{1-a}\quad &$a<1$,\\
  x_0 + Mb,\quad &$a=1$, \nonumber \\
 \label{limitingbehaviorend}
  a^M\left(x_0 + \frac{b}{a-1}\right) - \frac{b}{a-1},\quad &$a>1$.
\end{numcases}
For our three error sums, the constant $a$ is given by $a_2$ and $\tilde a_2$, and the constant $b$ is given by $b_2$, $b_3$ and $b_4$, respectively.
Applying the behavior of the iterated affine map $\psi_{a,b}$ as shown in \eqref{limitingbehaviorstart} to \eqref{limitingbehaviorend} to the maps \eqref{iteration1}, \eqref{iteration2} and \eqref{iteration3} with the constants $a_2$, $\tilde a_2$, $b_2$, $b_3$ and $b_4$, as well as the bounds \eqref{lemmabound2}, \eqref{lemmabound3} and \eqref{lemmabound4} shown in the three lemmas, we obtain the terms shown in \eqref{bound1}. \qed
\end{prf} 
\section{Numerical experiments.}\label{numericalexperiments}
For the case where the function $f$ associated to the operator $F$ in equation \ref{spde} is linear and given by $f(x)=x$, we can compute an exact solution by taking the function inside the linear operator $A$, so that we consider the eigenvalues of the operator $A+I$ which are given by $-\lambda_j+1$, $j\in\mathbb{N}$. For the linear function $f_1(x)=x$, we present convergence plots in Figure \ref{convergenceplotspeanut} and compare the approximated solution to the exact solution. For the nonlinear function $f_2(x) = \frac{1}{1+x^2}$, we additionally compare the approximate solution to a reference solution with $N=400$ and $M=100$ and show the convergence with respect to $N$ in the second plot. Since the constant for the noise term scales nonlinearly with the step size (see the final term in \eqref{finalscheme}), we cannot compute a reference solution to compare the other approximations to, as we cannot sensibly merge the finer random increments into more rough ones. For the convergence plots with $f = f_1$, the error was averaged over 100 independent realizations, while for $f=f_2$, 10 independent realizations were computed. We can see in all cases that the convergence with respect to $M$ and $N$ is of order one. We note at this point that we were not able to verify that the nonlinear function $f_2$ indeed satisfies the conditions \eqref{nonlinearbeginning} to \eqref{nonlinearend}, but convergence to a reference solution is observed nevertheless. We also remark that if strictly nonlinear $f$ is used, the computational effort is substantially higher since in every time step $N$ two-dimensional numerical integrals have to be computed.\\
 We now shed some light on how the error bound from Theorem \ref{strongerrortheorem} behaves here: We first note that for the Peanut shape it is $\lambda_1\approx 6.5155$, so in order for the condition $a_2<1$ to be fulfilled we would need $L\cdot N<6.5155.$ This is a very strict condition as either $N$ or $L$ needs to be quite small; for instance, if $N=100$, only nonlinearities with $L<0.065155$ are covered by the case $M\rightarrow \infty$ in Theorem \ref{strongerrortheorem}.  But as this is not a necessary condition to use Theorem \ref{strongerrortheorem}, we can still give bounds also if $a_2>1$.
\begin{example}\label{errorboundexample}
We give a concrete example for an error bound resulting from Theorem \ref{strongerrortheorem}. Suppose $T=0.1$, $N=100$, $M=50$, $R=301$, $L=0.01$, $\lambda_1 = 6.5155$, $|V_{k,j}^{(4)}|\leq 1$ for all $k=0,\dots, M$, $j=1,\dots, N$ (we have observed in practice that this has always been the case). We also assume $C=1$ (as Simpson's rule is very accurate for smooth functions) and $\epsilon_{k,j}^{(4)}\leq 1$ for all $k,j$ (which is a very loose assumption, but does not matter much as the corresponding terms in $b_3$ are small anyway). We assume error tolerances $\epsilon^{(\lambda)} = 2\cdot 10^{-4}$ and $\epsilon^{(\eta)} = 5\cdot 10^{-6}$. We also assume $\epsilon_0^{(2)} = \epsilon_0^{(3)} = \epsilon_0^{(4)} = 10^{-4}$ (as the initial condition is smooth, its approximation can be assumed to be accurate). With the above assumptions on the parameters, we obtain $a_2\approx 0.98904$ (with $\tilde a_2$, $b_2\approx 2.19316\cdot 10^{-6}$, $b_3 \approx 1.96257\cdot 10^{-4}$, $b_4 \approx 1.83156\cdot 10^{-4}$ and $b_5=5\cdot 10^{-5}$ and Theorem \ref{strongerrortheorem} yields a total error bound of $\mathcal{E}_M^{(234)}\leq 0.2234$. While this is clearly not a satisfactory bound especially compared to $\mathcal{E}_M^{(1)}$, we note that this bound depends critically on the error tolerances $\epsilon^{(\lambda)}$ and $\epsilon^{(\eta)}$. For example, for $\epsilon^{(\lambda)} = 2\cdot 10^{-5}$ and $\epsilon^{(\eta)} = 5\cdot 10^{-7}$, leaving all other parameters unchanged, we obtain $\mathcal{E}_M^{(234)}\leq 0.02257$. This shows that the error bound is usable in practice, but only with an as yet infeasible effort of computation (see Figure \ref{totalbasecomptime}).
\end{example}

In the plots in Figure \ref{exampleplots}, we give a visual impression of a solution of an SPDE \eqref{spde} on the Peanut shape by showing realizations of the approximated solution again for noise with $q_n = n^{-2}$, $A=\Delta$, $N=400$, $M=100$ on the asymmetric Peanut shape introduced in Section \ref{section_Numericalscheme}. For each simulation, the initial condition was picked to be the bump function
\begin{align*}
B_{x_1,x_2,y_1,y_2}(x,y) = \begin{cases}
  \exp\left( -\frac{1}{1-r^2}\right),\quad &r^2 = \left(x-\frac{x_1+x_2}{2}\right)^2 + \left(y-\frac{y_1+y_2}{2}\right)^2<1, \\
 0,\  &r^2 \geq 1
\end{cases}
\end{align*}
supported on an ellipse within the rectangle $[x_1,x_2]\times [y_1,y_2]$, where we picked $[x_1,x_2] = [0.4,0.6]$, $[y_1,y_2] = [0.3,0.5]$. The random input, i.e. the realizations of the Wiener increments was also the same for each simulation. Only the nonlinear function was changed for the different simulations: It was chosen to be a function
\begin{align}\label{bumpfunction}
    b_p(x) = \exp(-10\cdot(x-p)^2),
\end{align}
with one maximum at $p\in \{-0.4,-0.2,0.2,0.4\}$.
\begin{figure}
\newcommand{\Height}{5.3 cm}
\newcommand{\Width}{5.3 cm}
\newcommand{\X}{2.7}
\newcommand{\Yone}{1.0}
\newcommand{\Ytwo}{0.5}
\centering
 \begin{tabular}{@{}c@{}}
\begin{tikzpicture}
	\begin{axis}[
		height=\Height, 
		width=\Width,
		xmode = log,
		ymode = log,
		axis x line=bottom,
		axis y line=left,
		xlabel = {$M$},
		ylabel= {$\|  \tilde U_T^{M_{\mathrm{ref}},N} - \tilde U_T^{M,N}  \|$},
		clip=false,
		ymin=6.5e-6,ymax=4e-4,
		xmin=10,xmax=500,
		xtick={10,100},
		ytick={1e-4,1e-5},
		tick label style={font=\footnotesize},
		tick align=outside
		]
		\addplot[color=red,only marks,mark=x] table {Peanut_convergence_data_M.txt};
		\addplot[dashed,color=blue] table {referenceslope1.txt};
	\end{axis}
\end{tikzpicture}
\begin{tikzpicture}
	\begin{axis}[
		height=\Height, 
		width=\Width,
		ymode = log,
		xmode = log,
		axis x line=bottom,
		axis y line=left,
		xlabel = {$N$},
		ylabel= {$\|  \tilde U_T^{M,N_{\mathrm{ref}}} - \tilde U_T^{M,N}  \|$},
		clip=false,
		ymin=2e-4,ymax=1e-2,
		xmin=10,xmax=300,
		xtick={10,100},
		ytick={1e-2,1e-3,2e-4},
		yticklabels={{$1\cdot 10^{-2}$}, {$1\cdot 10^{-3}$},{$2\cdot 10^{-4}$}},
		tick label style={font=\footnotesize},
		y tick label style={/pgf/number format/sci},
		tick align=outside
		]
		\addplot[color=red,only marks,mark=x] table {Peanut_convergence_data_N.txt};
		\addplot[color=green,only marks,mark=o] table {Peanut_convergence_data_N_reference.txt};	
		\addplot[dashed,color=blue] table {referenceslope2.txt};
        \addlegendentry{$f = f_1$};
		\addlegendentry{$f = f_2$};
	\end{axis}
\end{tikzpicture}
\end{tabular}
\caption{Convergence tests for the Peanut shape. The first plot shows the convergence with respect to $M$ for the linear function $f_1(x)=x$, the second plot shows the convergence with respect to $N$ for $f_1(x)=x$ as well as $f_2(x) = \frac{1}{1+x^2}$. The dashed lines in blue are reference lines for an order of convergence equal to one.}
\label{convergenceplotspeanut}
\end{figure}
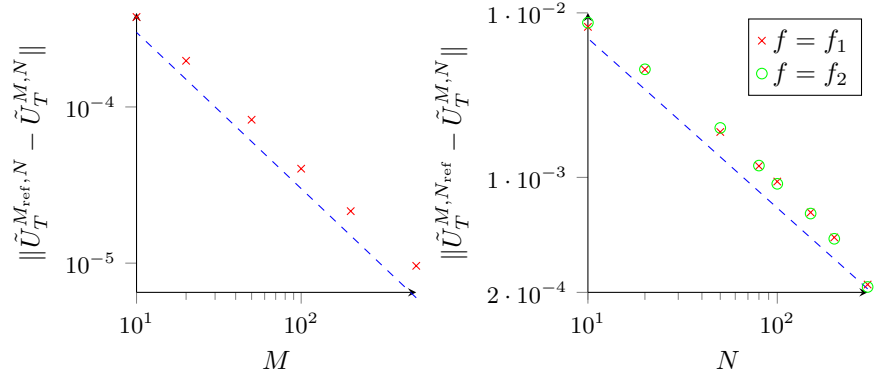
\begin{figure}[ht]
\newcommand{\Height}{5 cm}
\newcommand{\Width}{5 cm}
\newcommand{\X}{2.7}
\newcommand{\Yone}{1.0}
\newcommand{\Ytwo}{0.5}
\centering
 \begin{tabular}{@{}c@{}}
 \hspace*{-0.4em}{
 \begin{tikzpicture}
   \pgfdeclareimage[width=3.14cm]{img}{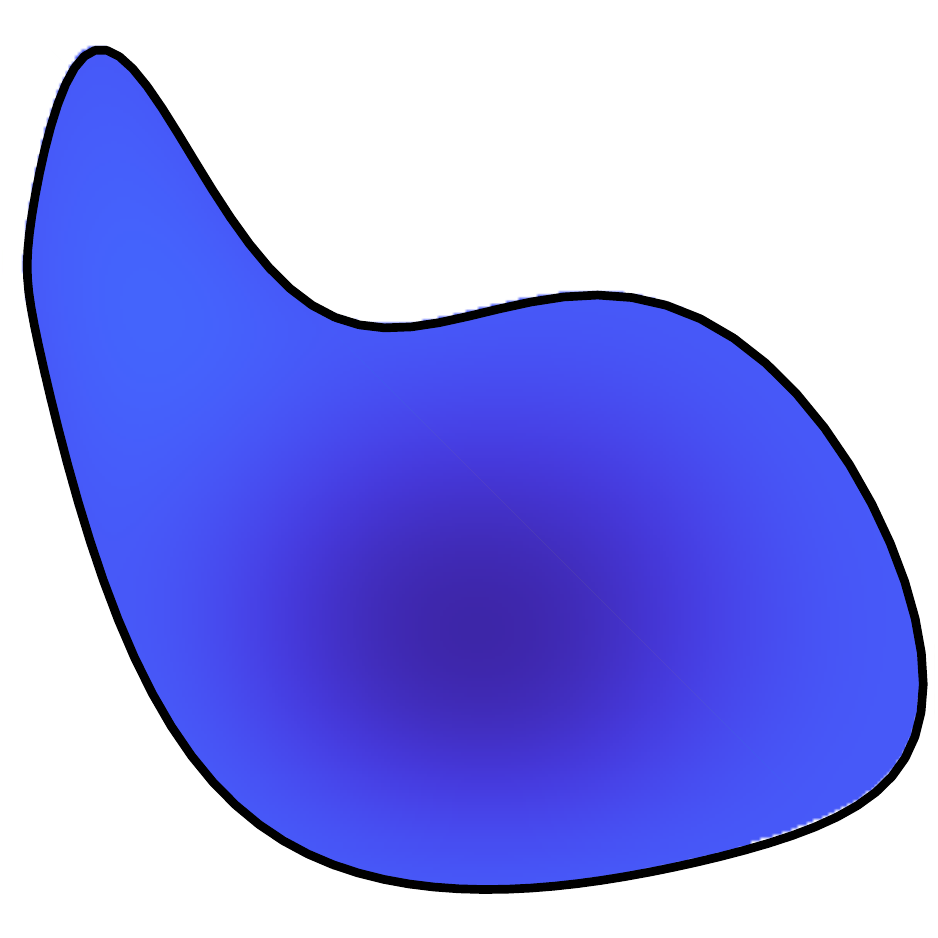}
        \node at (1.58,1.54) {\pgfuseimage{img}};
	\begin{axis}[
		height=\Height, 
		width=\Width,
		axis x line=bottom,
		axis y line=left,
		xlabel = {$x$},
		ylabel= {$y$},
		clip=false,
		ymin=0,ymax=1.1,
		xmin=0,xmax=1.1,
		xtick={0,0.2,0.4,0.6,0.8,1},
		ytick={0,0.2,0.4,0.6,0.8,1},
		tick label style={font=\footnotesize},
		tick align=outside
		]
	\end{axis}
\end{tikzpicture}
 }
 \begin{tikzpicture}
   \pgfdeclareimage[width=3.14cm]{img}{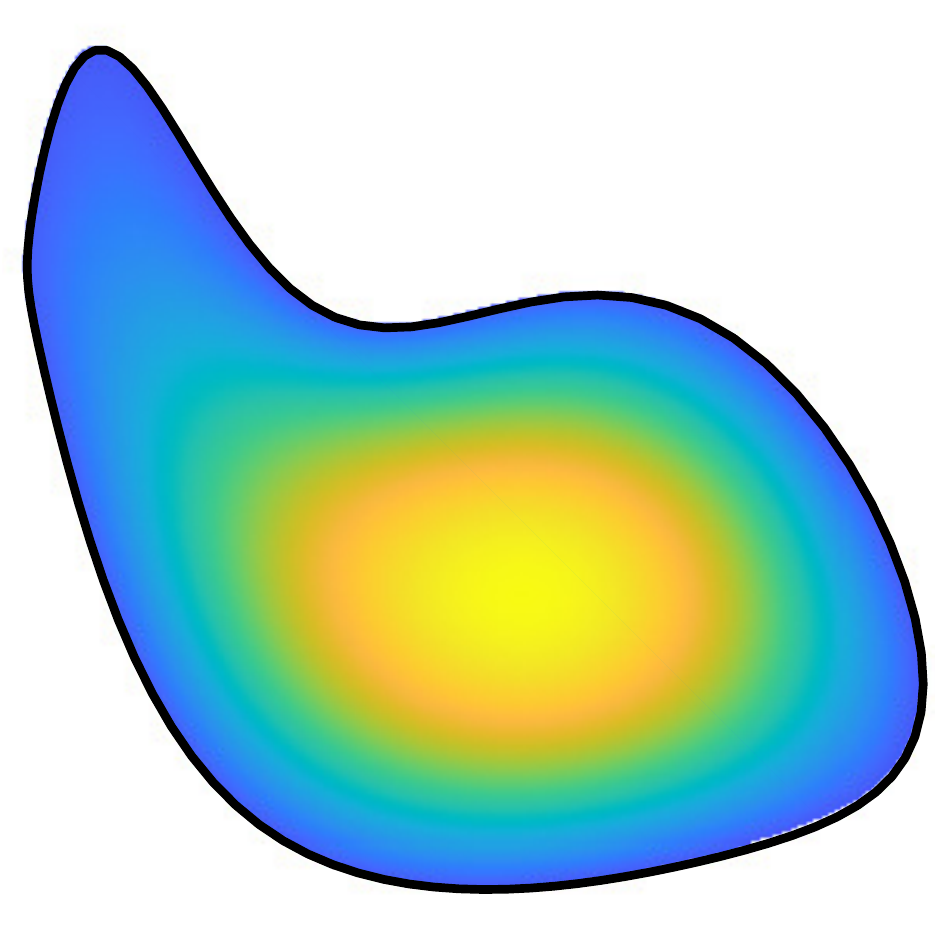}
        \node at (1.58,1.54) {\pgfuseimage{img}};
	\begin{axis}[
		height=\Height, 
		width=\Width,
		axis x line=bottom,
		axis y line=left,
		xlabel = {$x$},
		ylabel= {$y$},
		clip=false,
		ymin=0,ymax=1.1,
		xmin=0,xmax=1.1,
		xtick={0,0.2,0.4,0.6,0.8,1},
		ytick={0,0.2,0.4,0.6,0.8,1},
		tick label style={font=\footnotesize},
		tick align=outside
		]
	\end{axis}
\end{tikzpicture}\\
\hspace*{-0.7em}{
 \begin{tikzpicture}
   \pgfdeclareimage[width=3.14cm]{img}{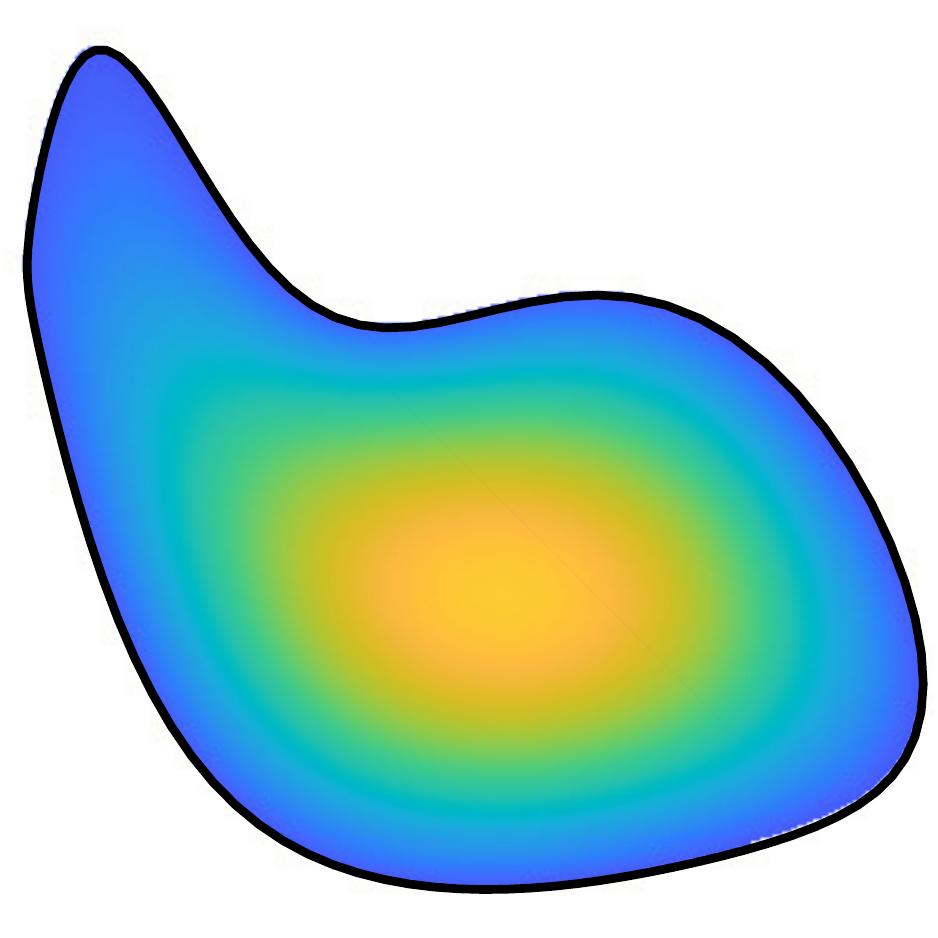}
        \node at (1.58,1.54) {\pgfuseimage{img}};
	\begin{axis}[
		height=\Height, 
		width=\Width,
		axis x line=bottom,
		axis y line=left,
		xlabel = {$x$},
		ylabel= {$y$},
		clip=false,
		ymin=0,ymax=1.1,
		xmin=0,xmax=1.1,
		xtick={0,0.2,0.4,0.6,0.8,1},
		ytick={0,0.2,0.4,0.6,0.8,1},
		tick label style={font=\footnotesize},
		tick align=outside
		]
	\end{axis}
\end{tikzpicture}
}
\begin{tikzpicture}
   \pgfdeclareimage[width=3.14cm]{img}{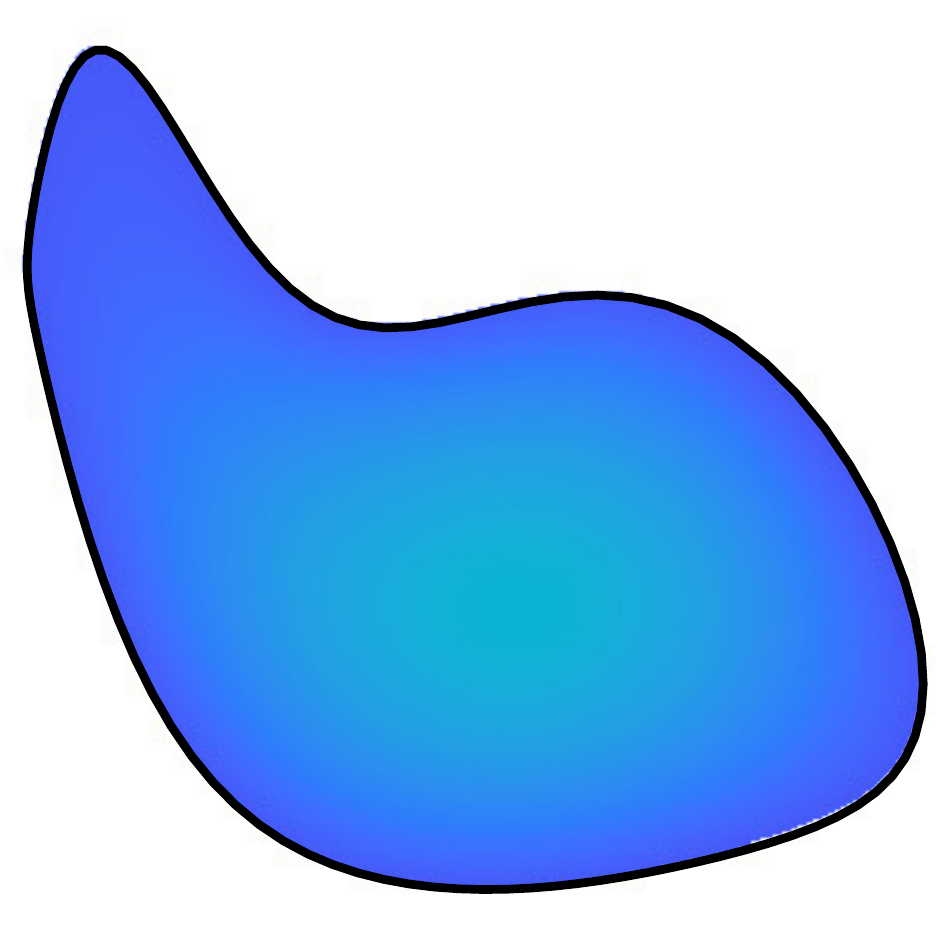}
        \node at (1.58,1.54) {\pgfuseimage{img}};
	\begin{axis}[
		height=\Height, 
		width=\Width,
		axis x line=bottom,
		axis y line=left,
		xlabel = {$x$},
		ylabel= {$y$},
		clip=false,
		ymin=0,ymax=1.1,
		xmin=0,xmax=1.1,
		xtick={0,0.2,0.4,0.6,0.8,1},
		ytick={0,0.2,0.4,0.6,0.8,1},
		tick label style={font=\footnotesize},
		tick align=outside
		]
	\end{axis}
\end{tikzpicture}
 \end{tabular}
 \centering
 \hspace*{1cm}
 \includegraphics[width=0.6\textwidth]{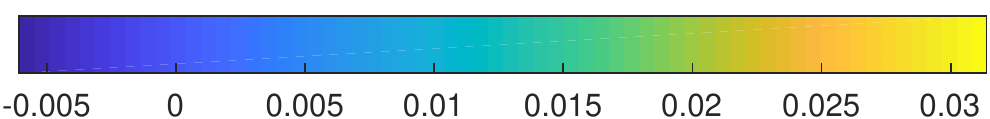}
\caption{The approximated solution of the SPDE \eqref{spde} on the Peanut shape for $N=400$, 
$T=0.1$, $M=100$ and the nonlinear functions \eqref{bumpfunction} for $p\in \{-0.4,-0.2,0.2,0.4\}$ from top-left to bottom-right. It appears that the process is `dimmed down' from using nonlinear functions with a peak of higher absolute value.}
\label{exampleplots}
\end{figure}
\section{Outlook.}\label{section_Outlook}
The advantages of the presented method are its flexibility in choosing the two-dimensional domain and the possibility to freely choose the points inside the domain at which the eigenfunctions are computed. This renders the solution of SPDEs with high resolution on a big class of domains possible. Due to the general setting of the spectral method, the scheme could be extended to a larger class of differential operators such as general second-order elliptic differential operators, whose Dirichlet eigenfunctions also form a complete orthogonal basis in $L^2(D)$ \cite{browder}. Difficulties might arise while using a boundary integral equation method, as explicit knowledge of the fundamental solution of the differential operator is required. The scheme could also be applied to three-dimensional domains with a $\mathcal{C}^2$ boundary: In order to do this, the fundamental solution needs to be changed to the one for the three-dimensional Helmholtz operator, the curved elements are now two-dimensional (see \cite[Section 5]{kleefeldlin}) and the linear system resulting from applying Beyn's contour integral algorithm is larger. However, for a small error, a very large number of eigenfunctions would have to be computed due to Weyl's law \eqref{weylslaw} and \eqref{errorjentzenpart}.\\
The requirements on the nonlinear operator $F$ could be relaxed by implementing a different scheme like the one from \cite{higherorder}, where the requirements on $F$ are less strict. Future work might also include deriving error bounds for higher moments such as $\mathbb{E}[\| U(t_k)-V_k^{(4)}\|^2]^{1/2}$.\\
%The acknowledgments section should not be numbered.
\section*{Acknowledgments}
We thank Arnulf Jentzen for valuable comments and a fruitful conversation.
\appendix
\section{Joint eigenbasis of $Q$ and $\Delta$ in two dimensions.}\label{appendixa}
Let $D\subset \mathbb{R}^d$ be a bounded domain, $d\geq 1$, with simple Dirichlet eigenvalues. It is well-known that a $Q$-Wiener process always has the series expansion given by \eqref{wienerseriesexpansion} in terms of the eigenfunctions and eigenvalues of $Q$ \cite[Proposition 4.3]{DaPrato}. In our setting (see end of Section \ref{section_Numericalscheme}), we have $\sum_{n=1}^\infty \|Q e^Q_n \|^2 = \sum_{n=1}^\infty q_n^2 = \sum_{n=1}^\infty n^{-(4+2\epsilon_q)}<\infty$, so $Q$ is a Hilbert-Schmidt operator. Furthermore, as is recalled in \cite[Section 3]{blomkernoise}, the operator $Q$ is Hilbert-Schmidt if and only if there is a kernel $q\in L^2(D\times D)$ such that
\begin{align}
    (Qf)(x) = \int_D q(x,y) f(y) \; \de y,\quad f\in L^2(D).
\end{align}
 It is shown in \cite[Theorem 3.3]{blomkernoise} that if $Q$ is trace class (in our setting, it is, since $\mathrm{tr}(Q) = \sum_{n=1}^\infty q_n<\infty$), the kernel $q$ can be understood as the correlation functional of a generalized Gaussian stochastic process (GSP) $\xi = \partial_t W$ satisfying
\begin{align}
    \mathbb{E}[\xi(t,x)] = 0,\quad \mathbb{E}[\xi(t,x)\xi(s,y)] = \delta(t-s)q(x,y),
\end{align}
with $0\leq s\leq t\leq T$ and $x,y\in D$. (We note here that this correspondence to a GSP relies on a geometric assumption \cite[Assumption 3.1]{blomkernoise} on $D$ which is not fulfilled by the Peanut shape considered in this paper, but it is also remarked on \cite[p. 259]{blomkernoise} that the results are true for `much more general domains'.) \\
A necessary and sufficient condition for $Q$ to have a joint eigenbasis with $\Delta$ is given in terms of $q$ \cite[Theorem 4.2]{blomkernoise}: If $q\in \mathcal{C}^2(D\times D)$, then $Q$ and $\Delta$ have a joint eigenbasis if and only if
\begin{align}\label{commutecondition1}
\Delta_y q(x,y) = \Delta_x q(x,y)\quad \text{and }\quad  q(\cdot,\zeta),\ q(\zeta,\cdot)\in \mathcal{C}_\tau^2(D)
\end{align}
for all $x,y,\zeta\in D$, where $\mathcal{C}^2_\tau (D) = \{ u\in \mathcal{C}^2(D)\ |\ u=0\ \text{on }\partial D\}$, and with $\Delta_x$ denoting the Laplacian in $x$ with Dirichlet boundary conditions and likewise for $y$.\\
Solving the first equation in \eqref{commutecondition1} by separation of variables, we write $q(x,y) = q_1(x)q_2(y)$ and obtain
\begin{align}
    q_1(x)\Delta_y q_2(y) = \Delta_x q_1(x) q_2(y) = \lambda \label{separated}
    \Leftrightarrow \frac{\Delta_y q_2(y)}{q_2(y)} = \frac{\Delta_x q_1(x)}{q_1(x)} = \lambda 
\end{align}
for a constant $\lambda$, for all $x,y\in D$ where $q_2(y)\neq 0$, $q_1(x)\neq 0$. We retrieve from \eqref{separated} the Helmholtz equations 
\begin{align}
    \Delta q_1(x) = \lambda q_1(x),\quad \Delta q_2(y) = \lambda q_2(y)
\end{align}
with homogeneous Dirichlet boundary conditions for $q_1$ and $q_2$, with the solutions $\{e_n\}_{n\in\mathbb{N}}$, $\{\lambda_n\}_{n\in\mathbb{N}}$ the Dirichlet eigenfunctions and eigenvalues on $D$, resulting in
\begin{align}
    q(x,y) = \left[\sum_{i=1}^\infty a_i e_i(x) \right]\cdot \left[\sum_{j=1}^\infty b_j e_j(y) \right] = \sum_{i=1}^\infty \sum_{j=1}^i a_i b_{i+1-j}e_i(x) e_{i+1-j}(y)
\end{align}
for real coefficients $a_n, b_n$, $n\in\mathbb{N}$. The first condition in \eqref{commutecondition1} now implies
\begin{align}
    \sum_{i=1}^\infty \sum_{j=1}^i \lambda_i a_i b_{i+1-j}e_i(x) e_{i+1-j}(y) = \sum_{i=1}^\infty \sum_{j=1}^i \lambda_{i+1-j} a_i b_{i+1-j}e_i(x) e_{i+1-j}(y)
\end{align}
and since the eigenvalues are simple, due to orthogonality only the terms remain where $i=i+1-j$. For real coefficients $c_n$, $n\in\mathbb{N}$, we obtain
\begin{align}\label{qseries}
    q(x,y) = \sum_{n=1}^\infty c_n e_n(x) e_n(y),
\end{align}
and in our setting in two dimensions we choose $c_n = n^{-(2+\epsilon_q)}$ (see Section \ref{setting}). Indeed, as $Q$ is symmetric and positive semidefinite \cite[p. 47 f.]{DaPrato}, and therefore so is $q$ \cite[p. 260]{blomkernoise}, due to Mercer's theorem \cite[p. 138]{couranthilbert}, $q$ always has a series representation \eqref{qseries} in terms of the eigenfunctions of $Q$, and we have shown that the Dirichlet eigenfunctions are the suitable choice for the Dirichlet boundary conditions.\\
On the second part in \eqref{commutecondition1}, we note that
\begin{align}\label{derivativeseries}
    \Delta_x q(x,y) = \sum_{n=1}^\infty n^{-(2+q_\epsilon)}\lambda_n e_n(x) e_n(y),
\end{align}
is absolutely convergent due to $n^{-(2+q_\epsilon)}\lambda_n = \mathcal{O}(n^{-(1+q_\epsilon)})$ being summable. We also note that the family $\{e_n\}_{n\in\mathbb{N}}$ is uniformly bounded in the $\|\cdot \|_\infty$ norm since due to the H\"{o}lder inequality
\begin{align}
    \|1\cdot e_n^2\|_1 \leq \| 1\|_1 \| e_n^2\|_\infty
    \Leftrightarrow \| e_n\|_2^2 \leq \mu(D)\cdot \| e_n\|_\infty^2,
\end{align}
so $\|e_n\|_\infty\leq \mu(D)^{-1/2}$, as $\|e_n\|_2=1$ for each $n\in\mathbb{N}$, where $\mu(D)$ is the area of $D$. It now follows from the Weierstra\ss\  convergence criterion \cite[p. 259]{forster} that \eqref{derivativeseries} converges absolutely and uniformly, and therefore as a uniform limit of continuous functions is continuous.

%%%%%%%%%%%%%%%%%%%%%%%%%%%%%%%%%%%%%%%%%%%%%%%%%%%%%%
%          7. REFERENCES SECTION
%%%%%%%%%%%%%%%%%%%%%%%%%%%%%%%%%%%%%%%%%%%%%%%%%%%%%%
%       READ THIS SECTION CAREFULLY
% Each of the references below MUST be cited in your article above. Do not include references that are not cited in your article.
% Follow the examples below carefully. We strongly suggest that you copy and paste your reference information directly into our examples.
% List all references in alphabetical order according to the first author's last name.
% Verify each URL works correctly and can be accessed properly. Your URL links should be to reputable websites. The command line for a website link begins with: \url{ }
% Do not add MR or DOI numbers to your references. AIMS production staff will add this information.
% Using BibTex is not recommended but can be handled.

\end{document}